\documentclass[notitlepage,leqno,10pt]{article}
\usepackage{graphicx,multirow}
\usepackage{amssymb,amsmath,amsthm}
\usepackage{lscape,lipsum,microtype}
\usepackage[margin=2.54cm]{geometry}  
\usepackage{hyperref}

\usepackage{ifluatex}
\ifluatex
 \usepackage{unicode-math}
\fi

\usepackage{latexsym}
\usepackage{amsmath}
\usepackage{amsfonts}
\usepackage{amssymb}
\usepackage{amscd}

\renewcommand{\a }{\alpha }

\renewcommand{\d}{\delta }
\newcommand{\D }{\Delta }

\renewcommand{\l }{\lambda }
\renewcommand{\L }{\Lambda }

\newcommand{\n }{\nabla }

\renewcommand{\O }{\Omega }

\newcommand{\ov}{\overline}
\newcommand{\intbar}{\mathop{\int\makebox(-13.5,0){\rule[4pt]{.7em}{0.3pt}}%
\kern-6pt}\nolimits}

\newcommand{\be}{\begin{equation}}
\newcommand{\ee}{\end{equation}}
\newcommand{\bes}{\begin{equation*}}
\newcommand{\ees}{\end{equation*}}
\newcommand{\ba}{\begin{eqnarray}}
\newcommand{\ea}{\end{eqnarray}}
\newcommand{\bas}{\begin{eqnarray*}}
\newcommand{\eas}{\end{eqnarray*}}
\newenvironment{pf}{\noindent{\sc Proof}.\enspace}{\rule{2mm}{2mm}\medskip}
\newenvironment{pfn}{\noindent{\sc Proof}}{\rule{2mm}{2mm}\medskip}

\newcommand{\R}{\mathbb{R}}

\newcommand{\Z}{\mathbb{Z}}

\newcommand{\N}{\mathbb{N}}

\author{Cheikh Birahim NDIAYE}
\date{}

\title{\bf Variational theory for the resonant \;$T$-curvature equation}

\begin{document}

\newtheorem{lem}{Lemma}[section]
\newtheorem{pro}[lem]{Proposition}
\newtheorem{thm}[lem]{Theorem}
\newtheorem{rem}[lem]{Remark}
\newtheorem{cor}[lem]{Corollary}
\newtheorem{df}[lem]{Definition}

\maketitle

\begin{center}
{\small

\noindent  Department of Mathematics Howard University \\  Annex 3, Graduate School of Arts and Sciences, \# 217 \\ DC 20059 Washington, USA

}

\end{center}

\footnotetext[1]{E-mail addresses: cheikh.ndiaye@howard.edu\\
\thanks{\\ The author was partially supported by NSF grant DMS--2000164.}}

\

\

\begin{center}
{\bf Abstract}
\end{center}
We study the resonant prescribed \;$T$-curvature problem on a compact \;$4$-dimensional Riemannian manifold with boundary. We derive sharp energy and gradient estimates of the associated Euler-Lagrange functional to characterize the critical points at infinity of the associated variational problem under a non-degeneracy on a naturally associated Hamiltonian function. Using this, we derive a Morse type lemma around the critical points at infinity. Using the Morse lemma at infinity, we prove new existence results of Morse theoretical type. Combining the Morse lemma at infinity and the Liouville version of the Barycenter technique of Bahri-Coron\cite{bc} developed in \cite{nd4}, we prove new existence results under a topological hypothesis on the boundary of the underlying manifold, and the entry and exit sets at infinity. 
 \begin{center}


\bigskip

\centerline{\bf AMS subject classification: 53C21, 35C60, 58J60, 55N10.}

\end{center}


\section{Introduction and statement of the results}

On a four-dimensional compact Riemannian manifolds with boundary \;$(\ov M, g)$, there exists a fourth-order operator \;$P_g$\; called Paneitz operator discovered by Paneitz\cite{p1} and an associated curvature quantity \;$Q_g$\; called  $Q$-curvature introduced by Branson-Oersted\cite{bo}. The Paneitz operator $P_g$ and the \;$Q$-curvature $Q_g$  are defined in terms of the Ricci tensor\;$Ric_{g}$\; and the scalar curvature\; $R_{g}$\; of  \;$(\ov{M}, g)$\; by
\begin{equation*}
P^4_g=\D_{g}^{2}-div_{g}\left((\frac{2}{3}R_{g}g-2Ric_{g})\n_g\right);\;\;\;\;\;\;\;Q_g=-\frac{1}{12}(\D_{g}R_{g}-R_{g}^{2}+3|Ric_{g}|^{2}),
\end{equation*}
where \;$div_g$\;is the divergence and\;$\n_g$\;is the covariant derivative of  with respect to \;$g$.
\vspace{4pt}

\noindent
On the other hand, Chang-Qing\cite{cq1} have discovered an operator\;$P^3_g$\; which is associated to the boundary \;$\partial M$\; of  \;$\ov M$\; and  a curvature quantity \;$T_g$\;naturally associated to \;$P^3_g$. They are defined by the formulas
\begin{equation*}
P^3_g=\frac{1}{2}\frac{\partial {\D_g}}{\partial n_g}+\D_{\hat g}\frac{\partial }{\partial n_g}-2H_g\D_{\hat g}+L_g(\nabla_{\hat g},\nabla_{\hat g})+\nabla_{\hat g}H_g.\nabla_{\hat g}+(F_g-\frac{R_g}{3})\frac{\partial }{\partial n_g}.
\end{equation*}
\begin{equation*}
T_g=-\frac{1}{12}\frac{\partial R_g}{\partial n_g}+\frac{1}{2}R_gH_g-<G_g,L_g>+3H_g^3-\frac{1}{3}tr_g(L_g^3)-\D_{\hat g}H_g,
\end{equation*}
where $\hat g$\;is the metric induced by \;$g$\;on\;\;$\partial M$,  $\frac{\partial }{\partial n_g}$ is the inward Neuman operator on $\partial M$ with respect to $g$,\;$L_g$ is the second fundamental form of \;$\partial M$ with respect to $g$,\;\;$H_g$\; is the mean curvature of \;$\partial M$\; with respect to \;$g$,\;$R^k_{g,ijl}$\;is the Riemann curvature tensor of \;$(\ov M, g)$,\;\; $R_{g,ijkl}=g_{mi}R^{m}_{g,jkl}$\;($g_{ij}$\;are the entries of the metric \;$g$),\;\;$F_g=R^{a}_{g,nan}$ (with $n$ denoting the index corresponding to the normal direction in local coordinates)\;\;\;and\;\;$<G_g,L_g>=\hat g^{ac}\hat g^{bd}R_{g,anbn}L_{g,cd}$. Moreover, the notation $L_g(\nabla_{\hat g},\nabla_{\hat g})$, means $L_g(\nabla_{\hat g},\nabla_{\hat g})(u)=\n_{\hat g}^a(L_{g,ab}\n_{\hat g}^bu)$. We point out that in all those notations above $i,j,k,l=1,\cdots 4$ and $a,b,c,d=1,\cdots 3$, and Einstein summation convention is used for repeated indices.
\vspace{6pt}

\noindent
As the Laplace-Beltrami operator and the Neumann operator are conformally covariant, we have that \;$P^4_g$ is conformally covariant of bidegree \;$(0, 4)$\ and \;$P^3_g$\; of bidegree \;$(0, 3)$. Furthermore, as they govern the transformation laws of the Gauss curvature and the geodesic curvature on compact surfaces with boundary, the couple \;$(P^4_g, P^3_g)$ does the same for \ $(Q_g,T_g)$\; on a compact four-dimensional Riemannian manifold with boundary \;$(\ov M, g)$. In fact, under a conformal change of metric\;$ g_u=e^{2u}g$, we have
\begin{equation}\label{eq:law}
\left\{
    \begin{split}
P^4_{g_u}=e^{-4u}P^4_g,\\
P^3_{g_u}=e^{-3u}P^3_{g},
    \end{split}
  \right.
\qquad \mbox{and}\qquad
\left\{
\begin{split}
P^4_gu+2Q_g=2Q_{ g_u}e^{4u}\;\;\text{in }\;\;M,\\
P^3_gu+T_g=T_{ g_u}e^{3u}\;\;\text{on}\;\;\partial M.
\end{split}
\right.
\end{equation}
Apart from this analogy, we have also an extension of the Gauss-Bonnet identity \;$\eqref{eq:gbc}$\; which is known as the Gauss-Bonnet-Chern formula
\begin{equation}\label{eq:gbc}
\int_{M}(Q_{g}+\frac{|W_{g}|^{2}}{8})dV_{g}+\oint_{\partial M}(T_g+Z_g)dS_g=4\pi^{2}\chi(\ov M)
\end{equation}
where\;$W_g$\;denote the Weyl tensor of\;$(\ov M,g)$\;and\;$Z_g$ is given by the following formula
$$
Z_g=R_gH_g-3H_gRic_{g,nn}+\hat g^{ac}\hat g^{bd}R_{g,anbn}L_{g,cd}-\hat g^{ac}\hat g^{bd}R_{g,acbc}L_{g,cd}+6H_g^3-3H_g|L_g|^2+tr_g(L_g^3),
$$
with \;$tr_g$ denoting the trace with respect to the metric induced on \;$\partial M$\; by \;$g$ (namely \;$\hat g$) and $\chi(\ov M)$ the Euler-Poincar\'e characteristic of \;$\ov M$. Concerning the quantity \;$Z_g$, we have that it vanishes when the boundary is totally geodesic and \;$\oint_{\partial M}Z_gdV_g$ is always conformally invariant, see \cite{cq1}. Thus, setting
\begin{equation}\label{eq:invc}
\kappa_{(P^4,P^3)}:=\kappa_{(P^4,P^3)}[g]:=\int_{M}Q_gdV_g+\oint_{\partial M}T_gdS_{g}
\end{equation}
we have that thanks to\;$\eqref{eq:gbc}$, and to the fact that \;$|W_g|^2dV_g$ is pointwise conformally invariant, $\kappa_{(P^4,P^3)}$\;is a conformal invariant (which justifies the notation used above). We remark that $4\pi^2$ is the the total integral of the $(Q, T)$-curvature of the standard four-dimensional Euclidean unit ball\; $\mathbb{B}^4.$
\vspace{8pt}

\noindent
As was asked in \cite{akn}, a natural question is whether every compact four-dimensional Riemannian manifold with  boundary\;$(\ov M, g)$\;carries a conformal metric \;$g_u$\; for which the corresponding \;$Q$-curvature \;$Q_{g_u}$\;is zero, the corresponding \;$T$-curvature \;$T_{g_u}$\; is a prescribed function and such that \;$(\ov M,g_u)$\;has minimal boundary. Thanks to $\eqref{eq:law}$, this problem is equivalent to finding a smooth solution to the following BVP:
\begin{equation*}
\left\{
\begin{split}
P^4_gu+2Q_g&=0\;\;&\text{in}\;\;M,\\
P^3_gu+T_g&=Ke^{3u}\;\;&\text{on}\;\;\partial M,\\
-\frac{\partial u}{\partial n_g}+H_gu&=0\;\;&\text{on}\;\;\partial M.
\end{split}
\right.
\end{equation*}
where \;$K: \partial M\longrightarrow \mathbb{R}_+$\; is a positive smooth function on \;$\partial M$.
\vspace{4pt}

\noindent
Since we are interested to find a metric in the conformal class of \;$g$, then we can assume\ that \;$H_g=0$, since this can be always obtained through a conformal transformation of the background metric. Thus, we are lead to solve the following BVP with Neumann homogeneous boundary condition:
\begin{equation}\label{eq:bvps}
\left\{
\begin{split}
P^4_gu+2Q_g&=0\;\;&\text{in}\;\;M,\\
P^3_gu+T_g&=Ke^{3u}\;\;&\text{on}\;\;\partial M,\\
\frac{\partial u}{\partial n_g}&=0\;\;&\text{on}\;\;\partial M.
\end{split}
\right.
\end{equation}
Defining \;$\mathcal{H}_{\frac{\partial}{\partial n}}$\;as
\begin{equation*}
\mathcal{H}_{\frac{\partial}{\partial n}}=\Big\{u\in W^{2, 2}(M):\;\;\;\frac{\partial u}{\partial n_g}=0\;\;\;\text{on}\;\;\partial M\Big\},
\end{equation*}
where \;$W^{2, 2}(M)$\; denotes the space of functions on \;$M$\; which are square integrable together with their first and second derivatives, and 
$$
\mathbb{P}^{4, 3}_g(u, v)=\left<P^4u, u\right>_{L^2(M)}+2\left<P^3_gu, u\right>_{L^2(\partial M)}, \;\;\;\;u,v\in\mathcal{H}_{\frac{\partial }{\partial n}},
$$
we have integration by part implies
\begin{equation}
\begin{split}
\mathbb{P}^{4, 3}_g(u, v)=\int_{M}\left(\D_g u\D_gv+\frac{2}{3}R_g\nabla_g u\cdot\nabla_g v\right)dV_g-2\int_{M}Ric_g(\nabla_g u,\nabla_g v)dV_g\\-2\oint_{\partial M}L_g( \nabla_{\hat g} u, \nabla_{\hat g} v)dS_g, \;\;\;u,v\in\mathcal{H}_{\frac{\partial }{\partial n}}
\end{split}
\end{equation}
 and is clearly a bilinear form on \;$\mathcal{H}_{\frac{\partial}{\partial n}}$. We set 
 \begin{equation}\label{kernl}
 \ker \mathbb{P}^{4, 3}_g:=\{u\in\mathcal{H}_{\frac{\partial }{\partial n}}:\;\;\;  \mathbb{P}^{4, 3}_g(u, v)=0, \;\;\forall v\in\mathcal{H}_{\frac{\partial }{\partial n}}\}
 \end{equation}
On the other hand, standard regularity theory implies that smooth solutions to \eqref{eq:bvps} can be found  by looking at critical points of the geometric functional
\begin{equation*}
\begin{split}
\mathcal{E}_g(u)=\mathbb{P}^{4, 3}_g(u, u)+4\int_{M}Q_gudV_g+4\oint_{\partial M}T_gudS_g-
\frac{4}{3}\kappa_{(P^4,P^3)}\log\oint_{\partial M}Ke^{3u}dS_g,\;\;\;u\in \mathcal{H}_{\frac{\partial }{\partial n}}.
\end{split}
\end{equation*}
As a Liouville type problem, the analytic features of equation \eqref{eq:bvps} and of the associated Euler-Lagrange functional \;$\mathcal{E}_g$\; depend strongly on the conformal invariant \;$\kappa_{(P^4,P^3)}$ . Indeed, depending on whether \;$\kappa_{(P^4,P^3)}$\; is a positive integer multiple of \;$4\pi^2$\; or not, the noncompactness of equation \eqref{eq:bvps} and the way of finding critical points of \;$\mathcal{E}_g$\; changes drastically. 
As far as existence questions are concerned, we have that problem \eqref{eq:bvps} has been solved in a work of Chang-Qing\cite{cq2} under the assumption that $\ker\; \mathbb{P}^{4, 3}_g\simeq\R$ , $\mathbb{P}^{4, 3}_g$\; is non-negative and \;$\kappa_{(P^4, \;P^3)}<4\pi^2$.  In \cite{nd6}, we show existence of solutions for \eqref{eq:bvps}  under the assumption \;$\ker \mathbb{P}^{4, 3}_g\simeq \R$ and \;$\kappa_{(P^4, \;P^3)}\notin 4\pi^2\N^*$.
\vspace{6pt}

\noindent
As a Liouviile type problem, the assumption \;$\ker \mathbb{P}^{4, 3}_g\simeq \R$ and \;$\kappa_{(P^4,P^3)}\notin 4\pi^2\N^*$ will be referred to as nonresonant case. This terminology is motivated by the fact that in that situation the set of solutions to some perturbations of equation (7) (including it) is compact. Naturally, we call resonant case when \;$\ker \mathbb{P}^{4, 3}_g\simeq \R$ and \;$\kappa_{(P^4, \;P^3)}\in 4\pi^2\N^*$. With these terminologies, we have that the works of Chang-Qing\cite{cq2} and our work in \cite{nd6} answer affirmatively the question raised above in the nonresonant case. However, for the resonant case, there are no known existence results to the best of our knowledge.
\vspace{10pt}

\noindent
In this work, beside existence results for \eqref{eq:bvps}, we are interested in a complete variational theory for  the boundary value problem \eqref{eq:bvps} in the resonant case, namely when \;$\ker \mathbb{P}^{4, 3}_g\simeq \R$ and \;$\kappa_{(P^4,P^3)}=4\pi^2k$\; for some \;$k\in \N^*$. To present the main results of the paper, we need to set first some notation and make some definitions. We define  the Hamiltonian function (at infinity) \;$\mathcal{F}_K: (\partial M)^k\setminus F_k(\partial M)\longrightarrow \R$\; by
$$
\mathcal{F}_K((a_1,\cdots, a_k)):=\sum_{i=1}^k\left(H(a_i, a_i)+\sum_{j=1,\; j\neq i} ^kG(a_i, \;a_j)+\frac{2}{3}\log(K(a_i))\right)
$$
where where \;$F_k(\partial M)$\; denotes the fat Diagonal of \;$(\partial M)^k$, namely $$F_k(\partial M):=\{A:=(a_1, \cdots, a_k)\in (\partial M)^k:\;\;\text{there exists} \;\;i\neq j\;\,\text{with}\;\, a_i=a_j\},$$ $G$ is the Green's function  defined by \eqref{green} , and $H$ is its regular defined as in \eqref{reg}. Furthermore, we define 
\begin{equation}\label{eq:critfk}
Crit(\mathcal{F}_K):=\{A\in (\partial M)^k\setminus F_k(\partial M), \;\;A\;\;\;\text{critical point of} \;\;\mathcal{F}_K\}.
\end{equation}
Moreover,  for $A=(a_1,\cdots, a_k)\in(\partial M)^k\setminus F_k(\partial M) $, and $i=1\cdots k$, we set
\begin{equation}\label{eq:partiallimit}
\mathcal{F}^A_i(x):=e^{3(H(a_i, x)+\sum_{j+1, \;j\neq i}^kG(a_j, x))+\frac{1}{3}\log(K(x))},
\end{equation}
and define
\begin{equation}\label{eq:defindexa}
\mathcal{L}_K(A):=-\sum_{i=1}^k  (\mathcal{F}^{A}_i)^{\frac{1}{2}}L_{\hat g}((\mathcal{F}^{A}_i)^{\frac{1}{6}})(a_i),
\end{equation}
where $$L_{\hat g}:=-\D_{\hat g}+\frac{1}{8}R_{\hat g}$$ is the conformal Laplacian associated to $\hat g$. We also set
\begin{equation}\label{eq:critsett}
\mathcal{F}_{\infty}:=\{A\in Crit(\mathcal{F}_K):\:\;\mathcal{L}_K(A)<0\},
\end{equation}
\begin{equation}\label{eq:minf}
i_{\infty}(A):=4k-1-Morse(A, \;\mathcal{F}_K),
\end{equation}
and define
\begin{equation}\label{eq:mi}
m_i^k:=card\{A\in Crit(\mathcal{F}_K):\; i_{\infty}(A)=i\}, \;\;i=0, \cdots, 4k-1,
\end{equation}
where \;$Morse(\mathcal{F}_K, A)$\; denotes the Morse index of \;$\mathcal{F}_K$\; at  \;$A$. We point out that for \;$k\geq 2$, $m_i^k=0$ for  $0\leq i\leq k-2$.\\
 For \;$k\geq 2$, we use the notation \;$B_{k-1}(\partial M)$\; to denote the set of formal barycenters of order \;$k-1$\; of \;$\partial M$, namely
\begin{equation}\label{eq:baryp}
B_{k-1}(\partial M):=\{\sum_{i=1}^{k-1}\alpha_i\d_{a_i}, \;\;a_i\in \partial M, \;\;\alpha_i\geq 0, i=1,\cdots, k-1,\;\,\sum_{i=1}^{k-1}\alpha_i=1\}.
\end{equation}
Furthermore, we define
\begin{equation}\label{eq:cpm}
c^{k-1}_p=dim \;H_p(B_{k-1}(\partial M)), \;\;p=1, \cdots 4k-5,
\end{equation}
where \;$H_p(B_{k-1}(\partial M)$\; denotes the \;$p$-th homology group of \;$B_{k-1}(\partial M)$\; with \;$\Z_2$ coefficients. Finally, we say
\begin{equation}\label{eq:nondeg}
(ND) \;\;\;\text{holds if } \;\;\mathcal{F}_K\;\text{is a Morse function and for every}\;\;A\in Crit(\mathcal{F}_K),  \;\;\mathcal{L}_K(A)\neq 0.
\end{equation}
\vspace{8pt}

\noindent
Now, we are ready to state our existence results  of Morse theoretical type  starting with the {\em critical} case, namely when \;$k=1$.
\begin{thm}\label{eq:morsepoincare1}
Let \;$(\ov M, \;g)$\; be a compact \;$4$-dimensional Riemannian manifold with \;boundary \;$\partial M$\; and interior \;$M$\; such that \;$H_g=0$, \;$\ker \mathbb{P}_g^{4, 3}\simeq \R$\; and \;$\kappa_{(P^4, \;P^3)}=4\pi^2$. Assuming that \;$K$\; is a smooth positive function on \;$\partial M$\; such that \;$(ND)$\, holds and the system
\begin{equation}\label{eq:mp1}
\begin{cases}
m^1_0=1+x_0,\\
m^1_i=x_i+x_{i-1}, \;&i=1, \cdots, 3,\\
0=x_3\\
x_i\geq 0,\;\;& i=0, \cdots, 3
\end{cases}
\end{equation}
has no solutions, then \;$K$\; is the \;$T$-curvature of a Riemannian metric  on $\ov M$\; conformally related to \;$g$ with zero \;$Q$-curvature in \;$M$\; and zero mean curvature on \;$\partial M$.
\end{thm}
\vspace{6pt}

\noindent
The system \eqref{eq:mp1} not having a solution traduces the violation of a strong Morse type inequalities ({\bf SMTI}) for the critical points at infinity of \;$\mathcal{E}_g$. Since ({\bf SMTI}) imply Poincare-Hopf type formulas,  then we have Theorem \ref{eq:morsepoincare1} implies the following Poincare-Hopf index type result.
\begin{cor}\label{eq:existence1}
Let \;$(\ov M, \;g)$\; be a compact \;$4$-dimensional Riemannian manifold with \;boundary \;$\partial M$\; and interior \;$M$\; such that \;$H_g=0$, \;$\ker \mathbb{P}_g^{4, 3}\simeq \R$\; and \;$\kappa_{(P^4, \;P^3)}=4\pi^2$.  Assuming that \;$K$\; is a smooth positive function on \;$\partial M$\; such that \;$(ND)$\, holds and
\begin{equation}\label{eq:ep1}
\sum_{A\in \mathcal{F}_{\infty}} (-1)^{i_{\infty}(A)}\neq 1,
\end{equation}
then \;$K$\; is the \;$T$-curvature of a Riemannian metric  on \;$\ov M$\; conformally related to \;$g$\; with zero \;$Q$-curvature in \;$M$\; and zero mean curvature on \;$\partial M$
 \end{cor}
\vspace{6pt}

\noindent
The formula \eqref{eq:ep1} says that the Euler number of the space of variations is different from the total contribution of the the true critical points at infinity and is of global character. Localizing the arguments of Corollary \ref{eq:existence1} in the case of the presence of a jump in the Morse index of the critical points of the Hamiltonian function \;$\mathcal{F}_K$, we have the following extension of Corollary \ref{eq:existence1}.
\begin{thm}\label{t:C}
Let \;$(\ov M, \;g)$\; be a compact \;$4$-dimensional Riemannian manifold with \;boundary \;$\partial M$\; and interior \;$M$\; such that \;$H_g=0$, \;$\ker \mathbb{P}_g^{4, 3}\simeq \R$\; and \;$\kappa_{(P^4, \;P^3)}=4\pi^2$ and \;$K$\; be a smooth positive function on \;$\partial M$\; satisfying the non degeneracy condition \;$(ND)$. Assuming that  there exists a positive integer \;$1 \leq l \leq 3$\; such that 
\begin{equation*}\label{eq:ep1c}
\begin{split}
\sum_{A\in \mathcal{F}_{\infty},\;  i_{\infty}(A) \leq l -1 } &(-1)^{i_{\infty}(A)}\neq 1\\&\text{and}\\
 \forall A \in \mathcal{F}_{\infty},\;\; &\quad i_{\infty}(A) \neq l, 
\end{split}
\end{equation*}
then \;$K$\; is the \;$T$-curvature of a Riemannian metric on \;$\ov M$\; conformally related to \;$g$ \;with zero \;$Q$-curvature in \;$M$\; and zero mean curvature on \;$\partial M$.

\end{thm}
\vspace{8pt}

\noindent
In the supercritical case, i.e \;$k\geq 2$, the Euler-Lagrange functional \;$\mathcal{E}_g$\; is not bounded from below, and taking into account the topological contribution of very large negative sublevels of \;$\mathcal{E}_g$, we have the following analogue of Theorem \ref{eq:morsepoincare1}.
\begin{thm}\label{eq:morsepoincare2}
Let \;$(\ov M, \;g)$\; be a compact \;$4$-dimensional Riemannian manifold with \;boundary \;$\partial M$\; and interior \;$M$\; such that \;$H_g=0$, \;$\ker \mathbb{P}_g^{4, 3}\simeq \R$, and \;$\kappa_{(P^4, P^3)}=4k\pi^2$\; with \;$k\geq 2$. Assuming that \;$K$\; is a smooth positive function on \;$\partial M$\; such that \;$(ND)$\, holds and the following system
\begin{equation}\label{eq:mp3}
\begin{cases}
0=x_0,\\
m^k_1=x_1,\\
m_i^k=c^{k-1}_{i-1}+x_i+k_{i-1}, \;&i=2,\cdots, 4k-4,\\
m_{i}^k=x_{i}+x_{i-1},\;& i=4k-3,\cdots, 4k-1,\\
0=x_{4k-1},\\
x_i\geq 0,\;\;& i=0, \cdots, 4k-1,
\end{cases}
\end{equation}
has no solutions,  then $\;K$\; is the \;$T$-curvature of a Riemannian metric on \;$\ov M$\; conformally related to \;$g$ \;with zero \;$Q$-curvature in \;$M$\; and zero mean curvature on \;$\partial M$.
\end{thm}
\begin{rem}\label{topnegc}
The presence of the number \;$c^{k-1}_{i-1}=dim H_{i-1}(B_{i-1}(\partial M))$\; in \eqref{eq:mp3} account for the contribution of the topology of very negative sublevels of \;$\mathcal{E}_g$. The relation between the topology of very negative sublevels of the Euler-Lagrange functional of Liouville type problems and the space of formal barycenters was first observed by Djadli-Malchiodi\cite{dm}.
\end{rem}
\vspace{6pt}

\noindent
As in the critical case, we have that Theorem \ref{eq:morsepoincare2} implies the following Poincar\'e-Hopf index type criterion for existence.
\begin{cor}\label{eq:existence2}
Let \;$(\ov M, \;g)$\; be a compact \;$4$-dimensional Riemannian manifold with \;boundary \;$\partial M$\; and interior \;$M$\; such that \;$H_g=0$, \;$\ker \mathbb{P}_g^{4, 3}\simeq \R$, and \;$\kappa_{(P^4, P^3)}=4k\pi^2$\; with \;$k\geq 2$. Assuming that \;$K$\; is a smooth positive function on \;$\partial M$\; such that \;$(ND)$\; holds and 
\begin{equation}\label{eq:ep2}
\sum_{A\in \mathcal{F}_{\infty}} (-1)^{i_{\infty}(A)}\neq \frac{1}{(k-1)!}\Pi_{i=1}^{k-1}(i-\chi(\partial M)),
\end{equation}
then \;$K$\; is the \;$T$-curvature of a Riemannian metric on \;$\ov M$\; conformally related to \;$g$ \;with zero \;$Q$-curvature in \;$M$\; and zero mean curvature on \;$\partial M$.

\end{cor}
\vspace{6pt}

\noindent
As in the critical case, we have that a localization of the arguments of Corollary \ref{eq:existence2} implies the following jumping index type result.
\begin{thm}\label{eq:Cm}
Let \;$(\ov M, \;g)$\; be a compact \;$4$-dimensional Riemannian manifold with \;boundary \;$\partial M$\; and interior \;$M$\; such that \;$H_g=0$, \;$\ker \mathbb{P}_g^{4, 3}\simeq \R$, \;$\kappa_{(P^4, P^3)}=4k\pi^2$\; with \;$k\geq 2$, and  let \;$K$\; be a smooth positive function on \;$\partial M$\; satisfying the non degeneracy condition  (ND).  Assuming that there exists a positive integer \;$1 \leq l \leq 4k-1$\;  and \;$A^l\in \mathcal{F}_{\infty}$ with \;$i_{\infty}(A^l) \leq l-1$ such that
\begin{equation*}\label{eq:ep2c}
\begin{split}
\sum_{A\in \mathcal{F}_{\infty}, \; i_{\infty}(A) \leq l - 1} (-1)^{i_{\infty}(A)}&\neq \frac{1}{(k-1)!}\Pi_{j=1}^{k-1}(j-\chi(\partial M))\\
&\text{and}\\
\forall A  \in \mathcal{F}_{\infty},\;\;  &\qquad i_{ \infty}(A)  \neq l,
\end{split}
\end{equation*} then \;$K$\; is the \;$T$-curvature of a Riemannian metric on \;$\ov M$\; conformally related to \;$g$ \;with zero \;$Q$-curvature in \;$M$\; and zero mean curvature on \;$\partial M$.

\end{thm}
\vspace{10pt}

\noindent
The Morse theoretical results stated above depend only the Morse Lemma at infinity around true critical points at infinity (see Lemma \ref{eq:morleminf}) which justify the condition \;$\mathcal{L}_K<0$\; in the definition of $\mathcal{F}_{\infty}$. However, our existence result of algebraic topological type are  based on the Morse lemma at infinity around all critical points at infinity. Thus, to state our existence result of algebraic topological type, we need first to introduce the neighborhood of potential critical points at infinity of \;$\mathcal{E}_g$. In order to do that, we first fix $\nu$ to be a positive and small real number, $\L$ to be a large positive constant, and \;$R$\; to be a large positive constant too. Next, for \;$\epsilon$\; small and positive, and \;$\Theta\geq 0$, we denote by \;$V(k, \epsilon,\Theta)$ the $(k,\epsilon, \Theta)$-{\em neighborhood of potential critical points at infinity}, namely
\begin{equation}\label{eq:lninfinity}
\begin{split}
&V(k, \epsilon, \Theta):=\{u\in\mathcal{H}_{\frac{\partial }{\partial n}}: \exists\; a_1, \cdots, a_k\in \partial M, \alpha_1,\cdots, \alpha_k>0,\;\;\l_1,\cdots, \l_k>0, \beta_1,\cdots, \beta_{\bar k}\in \R,\\& \;\;||u-\ov{u}_{(Q, T)}-\sum_{i=1}^k\alpha_i\varphi_{a_i,\lambda_i}-\sum_{r=1}^{\bar k}\beta_r(v_r-\ov{(v_r)}_{(Q, T)})||_{\mathbb{P}^{4, 3}}<\epsilon, \;\;\sum_{i=1}^k\alpha_i=k\;\;\;\alpha_i\geq 1-\nu,\\&\;\lambda_i\geq \frac{1}{\epsilon}, i=1,\cdots, k,\;\;\frac{2}{\Lambda}\leq \frac{\l_i}{\l_j}\leq \frac{\Lambda}{2},i, j=1, \cdots, k, \;|\beta_r|\leq \Theta, \;r=1, \cdots, \bar k, \\&\;\;\text{and}\;\;\l_i d_{\hat g}(a_i, a_j)\geq 4\ov CR\;\;\text{for}\;i\neq j\},
\end{split}
\end{equation}
where \;$\ov C$ is as in \eqref{eq:proua}, the \;$\varphi_{a_i, \l_i}$'s are as in \eqref{eq:projbubble}, \;$\bar k$\; is  as in \eqref{dimneg}, the \;$v_r$'s are  defined as in \eqref{eq:defeigen1},  the $\ov{(v_r)}_{(Q, t)}$'s are as in \eqref{eq:wmass}, and \;$||\cdot||_{\mathbb{P}^{4, 3}}$\; is defined as in \eqref{normp4p3}.
\vspace{6pt}

\noindent
As observed by Chen-Lin \cite{cl2} for Liouville type problems, the  minimization at infinity of Bahri-Coron \cite{bc} for Yamabe type problems has the following analogue for our problem. For \;$\Theta\geq 0$, there exists \;$\epsilon_0=\epsilon_0(\Theta)$ small and positive such that \;$\forall\;0<\epsilon\leq \epsilon_0$, we have
\begin{equation}\label{eq:mini}
\forall u\in V(k, \epsilon, \Theta), \text{the minimization problem }\\\min_{B_{\epsilon}^{\Theta}}||u-\ov{u}_{(Q, T)}-\sum_{i=1}^k\alpha_i\varphi_{a_i, \l_i}-\sum_{r=1}^{\bar k}\beta_r(v_r-\ov{(v_r)}_{(Q, T)})||_{\mathbb{P}^{4, 3}}
\end{equation}
has a unique solution, up to permutations, where \;$B_{\epsilon}^{\Theta}$\; is defined as follows
\begin{equation}
\begin{split}
B_{\epsilon}^{\Theta}:=&\{(\bar\alpha, A, \bar \l, \bar \beta)\in \R^k_+\times (\partial M)^k\times \R_+^k\times \R^{\bar k}:\sum_{i=1}^k\alpha_i=k,\; \a_i\geq 1-\nu,\;\l_i\geq \frac{1}{\epsilon}, i=1, \cdots, k,\\&\;|\beta_r|\leq \Theta, r=1, \cdots, \bar k, \;\l_id_{\hat g}(a_i, a_j)\geq 4\ov CR, i\neq j, i, j=1, \cdots, k\}.
\end{split}
\end{equation}
The selection map \;$s_k$\; is defined by \;$s_k: V(k, \epsilon, \Theta)\longrightarrow (\partial M)^k/\sigma_k$ as follows
\begin{equation}\label{eq:select}
s_k(u):=A, \;\;u\in V(k, \epsilon, \Theta), \;\,\text{and} \,\;A\;\;\text{is given by}\;\,\eqref{eq:mini},
\end{equation}
\vspace{8pt}

\noindent
We denote the critical points at infinity of \;$\mathcal{E}_g$\;  by \;$z^{\infty}$\; and use the notation \;$M_{\infty}(z^{\infty})$\; for their Morse indices at infinity, \;$W_u(z^{\infty})$\; for their unstable manifolds and \;$W_s(z^{\infty})$ for their stable manifolds, where $z$ is the corresponding critical point of $\mathcal{F}_K$. Furthermore, we denote by  \;$x^{\infty}$\; the "true" ones, namely \;$\mathcal{L}_K(x)<0$\; and the \;$y^{\infty}$\; the "false" ones, namely \;$\mathcal{L}_K(y)>0$. Moreover, we define \;$S$\; to be the following invariant set 
\begin{equation}
S:=\cup_{M_{\infty}(z_1^{\infty}), \;M_{\infty}(z_2^{\infty})\geq 4k-4+\bar k}\;W_u(z_1^{\infty})\cap W_s(z_2^{\infty}).
\end{equation}
We also define \;$S^{\infty}$\, to be the part of \;$S$\, at infinity, namely
\begin{equation}\label{eq:partinf}
S^{\infty}:=\cup_{M_{\infty}(z_1^{\infty}), \;M_{\infty}(z_2^{\infty})\geq 4k-4+\bar k}\;W_u^{\infty}(z_1^{\infty})\cap W_s(z_2^{\infty}),
\end{equation}
where \;$W_u^{–\infty}(z_1^{\infty})$\; denotes the restriction of \;$W_u(z_1^{\infty})$ at infinity. Furthermore, we denote by \;$S^{\infty}_{-}$\; the exit set from \;$S^{\infty}$\, starting from a false critical point at infinity \;$y^{\infty}$. \\
Similarly, we denote by \;$S^{\infty}_{+}$\; the entry set to \;$S^{\infty}$\; after having exited \;$S^{\infty}$\; through a set contained in \;$S^{\infty}_-$\; and entering into \;$S^{\infty}$\; through a true critical point at infinity \;$x^{\infty}$. 
Finally to state the result of algebraic topological type resulting from the Liouville version of the Barycenter technique of Bahri-Coron developed in \cite{nd4}, we need the existence of  \begin{equation}\label{eq:defom}
0\neq O^{*}_{\partial M}\in H^3(\partial M).
\end{equation}
Indeed, we prove:
\begin{thm}\label{eq:existence}
Let \;$(\ov M, \;g)$\; be a compact \;$4$-dimensional Riemannian manifold with \;boundary \;$\partial M$\; and interior \;$M$\; such that \;$H_g=0$, \;$\ker \mathbb{P}_g^{4, 3}\simeq \R$, and  \;$\kappa_{(P^4, P^3)}=4k\pi^2$\; with \;$k\geq 2$. Assuming that \;$K$\;is a smooth positive function on \;$\partial M$\; such that \;$(ND)$\; holds and either there is no \;$x^{\infty}$\; with \;$M_{\infty}(x^{\infty})=4k-4+\bar k$ \;or \;$s_k^*(O^*_{\partial M})\neq 0$\; in\; $H^3(S^{\infty})$\; and \;$s_k^*(O^*_M)=0$\; in\; $H^3(S_+^{\infty}\cup S_-^{\infty})$, then \;$K$\; is the \;$T$-curvature of a Riemannian metric conformally related to \;$g$ with zero \;$Q$-curvature in \;$M$\; and zero mean curvature on \;$\partial M$.
\end{thm}
\noindent
As in \cite{nd4}, Theorem \ref{eq:existence} implies the following collorary.
\begin{cor}\label{eq:cor1}
Let \;$(\ov M, \;g)$\; be a compact \;$4$-dimensional Riemannian manifold with \;boundary $\partial M$ and interior \;$M$\; such that \;$\ker \mathbb{P}_g^{4, 3}\simeq \R$\; and \;$\kappa_{(P^4, P^3)}=4k\pi^2$\; with \;$k\geq 2$. Assuming that $K$ is a smooth positive function on $\partial M$ such that \;$(ND)$\; holds and that every critical point $x$ of $\mathcal{F}_K$ of Morse index \;$0$\; or \;$1$\;satisfies\; $\mathcal{L}_K(x)<0$, then \;$K$\; is the \;$T$-curvature of a Riemannian metric conformally related to \;$g$ with zero \;$Q$-curvature in \;$M$\; and zero mean curvature \;$\partial M$.
\end{cor}.
\vspace{6pt}

\noindent
\section{Notation and preliminaries}\label{eq:notpre}
In the following, for a Riemmanian metric \;$\bar g$\; on \;$\partial M$ and $p\in \partial M$\;, we will use the notation \;$B^{\bar g}_{p}(r)$\; to denote the geodesic ball with respect to \;$\bar g$\; of radius \;$r$\;and center \;$p$. We also denote by \;$d_{\bar g}(x,y)$\; the geodesic distance with respect to \;$\bar g$\; between two points \;$x$\;and \;$y$\; of \;$\partial M$, $exp_x^{\bar g}$ the exponential map with respect to \;$\bar g$\; at \;$x\in \partial M$. $inj_{\bar g}(\partial M)$\;stands for the injectivity radius of \;$(\partial M, \bar g)$, $dV_{\bar g}$\;denotes the Riemannian measure associated to the metric\;$\bar g$.  Furthermore, we recall that \;$\n_{\bar g}$,  \;$\D_{\bar g}$, \;$R_{\bar g}$\;\; will denote respectively the covariant derivative, the Laplace-Beltrami operator, the scalar curvature and Ricci curvature with respect to \;$\bar g$. For simplicity, we will use the notation \;$B_p(r)$\; to denote $B^g_{p}(r)$, namely \;$B_p(r)=B^{\hat g}_p(r)$. $(\partial M)^2$\;stands for the cartesian product \;$\partial M\times \partial M$, while \;$Diag(\partial M)$\; is the diagonal of \;$(\partial  M)^2$.
\vspace{6pt}

\noindent
Similarly, for a Riemmanian metric \;$\tilde g$\; on \;$\ov M$\;, we will use the notation \;$B^{\tilde g, +}_{p}(r)$\; to denote the half geodesic ball with respect to \;$\tilde g$\; of radius \;$r$\;and center \;$p\in \partial M$. We also denote by \;$d_{\tilde g}(x,y)$\; the geodesic distance with respect to \;$\tilde g$\; between two points \;$x$\;and \;$y$\; of \;$\ov M$, $exp_x^{\tilde g}$ the exponential map with respect to \;$\tilde g$\; at \;$x\in \partial M$. $inj_{\tilde g}(\ov M)$\;stands for the injectivity radius of \;$(\ov M, g)$, $dV_{\tilde g}$\;denotes the Riemannian measure associated to the metric\;$g$, and \;$dS_{\tilde g}$\; the Riemannian measure associated to \;$\hat {\tilde g}:=\tilde g|_{\partial M}$, namely $dS_{\tilde g}=dV_{\hat {\tilde g}}$.  Furthermore, we recall that $\n_{\tilde g}$, \;$\D_{\tilde g}$, \;$R_{\tilde g}$\;\; will denote respectively the covariant derivative, the Laplace-Beltrami operator, the scalar curvature and Ricci curvature with respect to \;$g$. For simplicity, we will use the notation \;$B^+_p(r)$\; to denote $B^{g, +}_{p}(r)$, namely \;$B_p^+(r)=B^{g, +}_p(r)$, $p\in \partial M$
\vspace{6pt}

\noindent
For \;$1\leq p\leq \infty$\; and \;$k\in \N$, \;$\theta\in  ]0, 1[$, \;$L^p(M)$, $W^{k, p}(M)$, $C^k(\ov M)$, and $C^{k, \theta} (\ov M)$ stand respectively for the standard Lebesgue space, Sobolev space, $k$-continuously differentiable space and $k$-continuously differential space of H\"older exponent $\beta$, all with respect $g$.  Similarly, \;$1\leq p\leq \infty$\; and \;$k\in \N$, \;$\theta\in  ]0, 1[$, \;$L^p(\partial M)$, $W^{k, p}(\partial M)$, $C^k(\partial M)$, and $C^{k, \theta} (\partial M)$ stand respectively for the standard Lebesgue space, Sobolev space, $k$-continuously differentiable space and $k$-continuously differential space of H\"older exponent $\beta$, all with respect $\hat g$.
\vspace{4pt}

\noindent 
 Given a function \;$u\in L^1(M)\cap L^1(\partial  M)$, $\bar u_{\partial M}$\; and \;$\ov {u}_{(Q, T)}$\; are defined by  \;$$\bar u_{\partial M}=\frac{\oint_{\partial M} u(x)dS_{g}}{Vol_g(\partial M)},$$ with \;$Vol_{g}(\partial M)=\oint_{\partial M}dS_{g}$ and 
\begin{equation}\label{eq:wmass}
\ov{u}_{(Q, T)}=\frac{1}{4k\pi^2}\left(\int_{M} Q_gudV_{g}+\oint_{\partial M}T_gudS_g\right).
\end{equation}
\vspace{6pt}

\noindent
Given a generic Riemannian metric \;$\bar g$\; on \;$\partial M$\; and a function \;$F(x, y)$\; defined on a open subset of\;$(\partial M)^2$\; which is symmetric and\; with $F(\cdot, \cdot)\in C^2$ with respect to \;$\bar g$, we define \;$\frac{\partial F(a, a)}{\partial a}:=\frac{\partial F(x, a)}{\partial x}|_{x=a}=\frac{\partial F(a, y)}{\partial y}|_{y=a}=$, and \;$\D_{\bar g} F(a_1, a_2):=\D_{\bar g, x}F(x, a_2)|_{x=a_1}=\D_{\bar g, y}F(a_2, y)|_{y=a_1}.$
\vspace{6pt}

 \noindent
 For  \;$\epsilon>0$\; and small,  $\l\in \R_+$, $\l\geq \frac{1}{\epsilon}$, and \;$a\in \partial M$, $O_{\l, \epsilon}(1)$\; stands for quantities bounded uniformly in \;$\l$, and $\epsilon$, and \;$O_{a, \epsilon}(1)$ stands for quantities bounded uniformly in $a$ and $\epsilon$. For  $l\in \N^*$, $O_{l}(1)$ stands for quantities bounded uniformly in \;$l$\; and \;$o_l(1)$ stands for quantities which tends to $0$ as $l\rightarrow +\infty$.  For  $\epsilon$ positive and small, \;$a\in \partial M$\; and \;$\l\in \R_+$ large, $\l\geq \frac{1}{\epsilon}$,\;$O_{a, \l, \epsilon}(1)$\; stands for quantities bounded uniformly in \;$a$, \;$\l$, and $\epsilon$. For $\epsilon$ positive and small, $p\in \N^{*}$, $\bar \l:=(\l_1, \cdots, \l_p)\in (\R_+)^p$, $\l_i\geq \frac{1}{\epsilon}$  for $i=1, \cdots, p$, and $A:=(a_1, \cdots, a_p)\in M^p$ (where $ (\R_+)^p$ and \;$(\partial M)^p$\; denotes respectively the cartesian product of $p$ copies of $\R_+$ and $\partial M$), $O_{A, \bar \l, \epsilon}(1)$ stands for quantities bounded uniformly in $A$, $\bar \l$, and $\epsilon$. Similarly for $\epsilon $ positive and small,  $p\in \N^{*}$, $\bar \l:=(\l_1, \cdots, \l_p)\in (\R_+)^p$, $\l_i\geq \frac{1}{\epsilon}$ for $i=1, \cdots, p$, $\bar \alpha:=(\alpha_1, \cdots, \alpha_p)\in \R^p$, $\alpha_i$ close to $1$ for $i=1, \cdots, p$, and $A:=(a_1, \cdots, a_p)\in (\partial M)^p$ (where $ \R^p$  denotes the cartesian product of $p$ copies of $\R$, $O_{\bar\alpha, A, \bar \l, \epsilon}(1)$ will mean quantities bounded from above and below independent of $\bar \alpha$, $A$, $\bar \l$, and $\epsilon$. For $x\in \R$, we will use the notation $O(x)$ to mean $|x|O(1)$ where $O(1)$ will be specified in all the contexts where it is used. Large positive constants are usually denoted by $C$ and the value of\;$C$\;is allowed to vary from formula to formula and also within the same line. Similarly small positive constants are also denoted by $c$ and their value may varies from formula to formula and also within the same line.
\vspace{6pt}

\noindent
We say \;$\mu\in \R$\; is an eigenvalue of the \;$P_g^4$\; to \;$P^3_g$\; operator on \;$\mathcal{H}_{\frac{\partial }{\partial n}}$\; if there exists \;$0\neq v\in  W^{2, 2}(M)$\; such that
\begin{equation}\label{eq:defeigen0}
\left\{
\begin{split}
P^4_gv&=0\;\;&\text{in}\;\;M,\\
P^3_gv&=\mu v\;\;&\text{on}\;\;\partial M,\\
\frac{\partial v}{\partial n_g}&=0\;\;&\text{on}\;\;\partial M.
\end{split}
\right.
\end{equation}
By abuse of notation, we call \;$v$\; in \eqref{eq:defeigen0} an eigenfunction associated to \;$\mu$. We call \;$\bar k$\; the number of negative eigenvalues (counted with multiplicity) of  the \;$P_g^4$\; to \;$P^3_g$ operator on \;$\mathcal{H}_{\frac{\partial }{\partial n}}$. We point out that \;$\bar k$\; can be zero, but it is always finite. If \;$\bar k\geq 1$, then we will denote by \;$E_{-}\subset \mathcal{H}_{\frac{\partial }{\partial n}}$ the direct sum of the eigenspaces corresponding to the negative eigenvalues of  the \;$P_g^4$\; to \;$P^3_g$\; operator\; on \;$\mathcal{H}_{\frac{\partial }{\partial n}}$ . The dimension of \;$E_{-}$ is of course $\bar k$, i.e
\begin{equation}\label{dimneg}
\bar k=\dim E_{-}.
\end{equation}
On the other hand, we have the existence of a basis of eigenfunctions \;$v_1,\cdots, v_{\bar k}$\; of \;$E_{-}$ satisfying
\begin{equation}\label{eq:defeigen1}
\left\{
\begin{split}
P^4_gv_r&=0\;\;&\text{in}\;\;M,\\
P^3_gv_r&=\mu_r v_r\;\;&\text{on}\;\;\partial M,\\
\frac{\partial v_r}{\partial n_g}&=0\;\;&\text{on}\;\;\partial M.
\end{split}
\right.
\end{equation}
\begin{equation}\label{eq:defeigen2}
\mu_1\leq \mu_2\leq \cdots\leq \mu_{\bar k}<0<\mu_{\bar k+1}\leq\cdots,
\end{equation}
where \;$\mu_r$'s are the eigenvalues of  the operator \;$P_g^4$\; to \;$P^3_g$\; on \;$\mathcal{H}_{\frac{\partial }{\partial n}}$ counted with multiplicity. 
We define \;$\mathbb{P}^{4, 3}_{g, +}$\; by
\begin{equation}\label{eq:nonnegative}
\mathbb{P}^{4, 3}_{g, +}(u, v)=\mathbb{P}^{4, 3}_{g}(u, v)-2\sum_{r=1}^{\bar k}\mu_r\left(\oint_{\partial M} uv_rdS_g\right) \left(\oint_{\partial M} vv_rdS_g\right).
\end{equation}
$\mathbb{P}^{4, 3}_{g, +}$\; is obtained by just reversing the sign of the negative eigenvalue of \;$\mathbb{P}^{4, 3}_{g}$. 
We set also
\begin{equation}\label{normp4p3}
||u||_{\mathbb{P}^{4, 3}}:=\sqrt{\mathbb{P}^{4, 3}_{g, +}(u, u)}, \;\;\;\text{and}\;\;\;\left<u, v\right>_{\mathbb{P}^{4, 3}}=P^{4, 3}_{g, +}(u, v),
\end{equation}
where \;$\mathbb{P}^{4, 3}_{g, +}$\; is defined as in \eqref{eq:nonnegative}. We can choose \;$v_1,\cdots, v_{\bar k}$\; so that they constitute a \;$\left<\cdot, \cdot\right>_{\mathbb{P}^{4, 3}}$-orthonormal basis fo \;$E_-$. We denote by \;$\n^{\mathbb{P}^{4, 3}}$\; the gradient with respect to $\left<\cdot, \cdot\right>_{\mathbb{P}^{4, 3}}$.
\vspace{6pt}

\noindent
For \;$t>0$, we define the following perturbed functional
\begin{equation}\label{eq:jt}
(\mathcal{E}_g)_t(u):=\mathbb{P}^{4,3}(u, u)+4t\int_{M}Q_gudV_g+4t\oint_{\partial M}T_gudS_g-
\frac{4}{3}t\kappa_{(P^4,P^3)}\log\oint_{M}Ke^{3u}dS_g,\;\;u\in \mathcal{H}_{\frac{\partial }{\partial n}}.
\end{equation}
\vspace{6pt}

\noindent
$\bar B^{\bar k}_r$ will stand for the closed ball of center \;$0$ and radius \;$r$\; in \;$\R^{\bar k} $. $\mathbb{S}^{\bar k-1}$\; will denote the boundary of \;$\bar B^{\bar k}_1$. Given a set \;$X$,  we define \;$\widetilde{X\times \bar B^{\bar k}_1}$\; to be the cartesian product \;$X\times \bar B^{\bar k}_1$\; where the tilde means that \;$X\times \partial B^{\bar k}_1$\; is identified with \;$\partial B_1^{\bar k}$.
\vspace{6pt}

\noindent
In the sequel also, \;$(\mathcal{E}_g)^{c}$\; with \;$c\in \R$ will stand for \;$(\mathcal{E}_g)^{c} :=\{u\in \mathcal{H}_{\frac{\partial }{\partial n}}:\;\;\mathcal{E}_g(u)\leq c\}$. For \;$X$\; a topological space, $H_{*}(X)$\; will denote the singular homology of \;$X$, $H^*(X)$ for the cohomology, and $\chi(X)$\; the Euler characteristic of \;$X$, all with \;$\Z_2$\; coefficients.

\vspace{4pt}

\noindent
As above, in the general case, namely \;$\bar k\geq 0$, for $\epsilon$\; small and positive, $\bar \beta:=(\beta_1, \cdots, \beta_{\bar k})\in \R^{\bar k}$\;  with \;$\beta_i$\; close to \;$0$, $i=1, \cdots, \bar k$) (where \;$\R^{\bar k}$\; is the empty set when \;$\bar k=0$), $\bar \l:=(\l_1, \cdots, \l_p)\in (\R_+)^p$, $\l_i\geq \frac{1}{\epsilon}$\;  for $i=1, \cdots, p$, $\bar \alpha:=(\alpha_1, \cdots, \alpha_p)\in \R^p$, $\alpha_i$ close to $1$ for $i=1, \cdots, p$, and $A:=(a_1, \cdots, a_p)\in (\partial M)^p$, $p\in \N^{*}$, $w\in \mathcal{H}_{\frac{\partial }{\partial n}}$\; with \;$||w||_{\mathbb{P}^{4, 3}}$\; small, $O_{\bar\alpha, A, \bar \l, \bar \beta, \epsilon}(1)$\; will stand quantities bounded independent of \;$\bar \alpha$, $A$, $\bar \l$, $\bar \beta$, and \;$\epsilon$, and $O_{\bar\alpha, A, \bar \l, \bar \beta, w, \epsilon}(1)$\; will stand quantities bounded independent of \;$\bar \alpha$, $A$, $\bar \l$, $\bar \beta$,  $w$\; and \;$\epsilon$.
\vspace{6pt}

\noindent
For point \;$b\in \R^3$\; and \;$\lambda$\; a positive real number, we define \;$\delta_{b, \lambda}$\; by \begin{equation}\label{eq:standarbubble}
\delta_{b, \lambda}(y):=\log\left(\frac{2\lambda}{1+\lambda^2|y-b|^2}\right),\;\;\;\;\;\;y\in \R^3.
\end{equation}
The functions \;$\delta_{b, \lambda}$\; verify the following equation
\begin{equation}\label{eq:bubbleequation}
(-\D_{\R^3})^{\frac{3}{2}}\delta_{b,\lambda}=2e^{3\delta_{b,\lambda}}\;\;\;\text{in}\;\;\;\R^3.
\end{equation}
Using the existence of conformal Fermi coordinates, we have that, for \;$a \in \partial M$\; there exists a function \;$u_a\in C^{\infty}(\ov M)$ such that
\begin{equation}\label{eq:detga}
g_a = e^{2u_a} g\;\; \text{verifies}\;\;det g_a(x)=1+O(d_{g_a}(x, a)^m)\;\;\text{for}\;\;\; x\in B^{g_a, +}_a( \varrho_a).
\end{equation}
with \;$0<\varrho_a<\frac{inj_{g_a}(\ov M)}{10}$. Moreover, we can take the families of functions \;$u_a$, $g_a$ and $\varrho_a$ such that
\begin{equation}\label{eq:varro0}
\text{the maps}\;\;\;a\longrightarrow u_a, \;g_a\;\;\text{are}\;\;C^1\;\;\;\text{and}\;\;\;\;\varrho_a\geq \varrho_0>0,
\end{equation}
for some small positive \;$\varrho_0$\; satisfying \;$\varrho_0<\inf\{\frac{inj_g(M)}{10}, \frac{inj_{\hat g}(\partial M)}{10}\}$, and
\begin{equation}\label{eq:proua}
\begin{split}
&||u_a||_{C^4(\ov M)}=O_a(1),\;\;\frac{1}{\ov C^2} g\leq g_a\leq \ov C^2 g, \\\;&u_a(x)= O_a(d^2_{\hat g_a}(a, x))=O_a(d_{\hat g}^2(a, x)) \;\;\text{for}\;\;x\in\;\;B_a^{\hat g_a}(\varrho_0)\supset B_a(\frac{\varrho_0}{2\ov C}),\;\;\text{and}\\&
u_a(a)=0,\;\;\;R_{\hat g_a}(a)=0, \;\;\;\frac{\partial u_a}{\partial n_g}(a)=0,
\end{split}
\end{equation}
for some large positive constant \;$\ov C$\; independent of \;$a$. 
For \;$a\in \partial M$, and \;$r>0$, we set
\begin{equation}\label{eq:expballaga}
exp_a^{a}:=exp_{a}^{\hat g_a}\;\;\;\;\text{and}\;\;\;B_a^{a}(r):=B_a^{\hat g_a}(r).
\end{equation}
Now, for \;$0<\varrho<\min\{\frac{inj_g(\partial M)}{4}, \frac{\varrho_0}{4}\}$\; where \;$\varrho_0$\; is as in \eqref{eq:varro0}, we define a smooth cut-off function  satisfying the following properties:
\begin{equation}\label{eq:cutoff}
\begin{cases}
\chi_{\varrho}(t)  = t \;\;&\text{ for } \;\;t \in [0,\varrho],\\
\chi_{\varrho}(t) = 2 \varrho \;\;&\text{ for } \;\; t \geq 2 \varrho, \\
 \chi_{\varrho}(t) \in [\varrho, 2 \varrho] \;\;\;&\text{ for } \;\; t\in [\varrho, 2 \varrho].
\end{cases}
\end{equation}
Using the cut-off function $\chi_{\varrho}$, we define for\; $a\in \partial M$\; and \;$\lambda\in \R_+$\; the function \;$\hat{\delta}_{a, \lambda}$\; as follows
\begin{equation}\label{eq:hatdelta}
\hat{\delta}_{a, \lambda}(x):=\log\left(\frac{2\lambda}{1+\lambda^2\chi_{\varrho}^2(d_{\hat g_a}(x, a))}\right).
\end{equation}
For every \;$a\in \partial M$\; and \;$\lambda\in \R_+$, we define \;$\varphi_{a, \lambda}$\; to be the unique the solution of
\begin{equation}\label{eq:projbubble}
\begin{cases}
P^4_g \varphi_{a, \l}+\frac{2}{k}Q_g=0\;\;\;\text{in}\;\;\;M,\\
P_g^3 \varphi_{a, \l} \, + \, \frac{1}{k} T_g\, = 4\pi^2\,\frac{ e^{3 (\hat{\d}_{a,\l} \,  +  \, u_a)}}{\oint_M e^{3 (\hat{\d}_{a,\l} \,  +  \, u_a)}dS_g}\;\; \mbox{ in } \;\;\partial M, \\
\frac{\partial \varphi_{a, \l}}{n}=0,\\
(\phi_{a, \l})_{(Q, T)}=0.
\end{cases}
\end{equation}
Next,  let \;$S (a, x)$\,; $(a, x)\in \partial M\times \ov M$\; be defined by 
\begin{equation}\label{eq:defG4}
\begin{cases}
P_g^4 S(a, \cdot)+\frac{2}{k}Q_g(\cdot)=0\;\;\text{in}\;\; M,\\
P^3_gS(a, \cdot)+\frac{1}{k}T_g(\cdot)= 4\pi^2\d_a(\cdot),\;\;\text{on}\;\; \partial M,\\
\frac{\partial S (a, \cdot)}{\partial n} =0\;\;\;\text{on}\;\; \partial M,\\
\int _MS(a, x)Q_g(x)dV_g(x)=0.
\end{cases}
\end{equation}
Then
 \begin{equation}\label{green}
G(a, \cdot)=S(a, \cdot)|_{\partial M}.
\end{equation} 
is a Green's function of the \;$P^4_g+\frac{2}{k}Q_g(\cdot)$\; to \;$P^3_g+\frac{1}{k}T_g(\cdot)$\; operator on \;$\mathcal{H}_{\frac{\partial }{\partial n}}$.
Thus, we have the integral representation: $\forall u\in \mathcal{H}_{\frac{\partial }{\partial n}}$\; such that \;$P^4_g u+\frac{2}{k}Q_g=0$, 
\begin{equation}\label{eq:G4integral}
u(x)-\ov{u}_{(Q, T)}=\frac{1}{4\pi^2}\oint_{\partial M} G(x, y)P_g^3u(y), \;\; \;x\in \partial M.
\end{equation}
Moreover, \;$G$\; decomposes as follows
\begin{equation}\label{reg}
 G(a,x)=\log \left(\frac{1}{\chi_{\varrho}^2(d_{{\hat g}_a}(a, x))}\right)+H(a, x),
\end{equation}
where \;$H$\; is the regular par of \;$G$. Furthermore, we have
\begin{equation}\label{eq:regH4}
G\in C^{\infty}((\partial M)^2-Diag(\partial M)),\;\;\;\;\text{and} \;\;\;H\in C^{3, \beta}((\partial M)^2)\;\;\;\forall \beta\in (0, 1).
\end{equation}
B symmetry of \;$H$, we have
\begin{equation}\label{eq:relationderivative}
\frac{\partial \mathcal{F}(a_1, \cdots, a_k)}{\partial a_i}=\frac{2}{3}\frac{\n_{\hat g}\mathcal{F}^{A}_i(a_i)}{\mathcal{F}^{A}_i(a_i)}, \;\;\;i=1, \cdots, k.
\end{equation}
Next, setting
\begin{equation}\label{eq:deflA}
l_K(A):=\sum_{i=1}^k\left(\frac{\D_{\hat g} \mathcal{F}^{A}_i(a_i)}{(\mathcal{F}^{A}_i(a_i))^{\frac{1}{3}}}-\frac{3}{4}R_{\hat g}(a_i)(\mathcal{F}^{A}_i(a_i))^{\frac{2}{3}}\right),
\end{equation}
we have 
\begin{equation}\label{eq:auxiindexa1}
l_K(A)=6\mathcal{L}_K(A), \;\;\forall A\in Crit(\mathcal{F}_K).
\end{equation}
For \;$k\geq 2$, we denote by \;$B_k(\partial M)$\; the set of formal barycenters of \;$\partial M$\; of order \;$k$, namely
 \begin{equation}\label{eq:barytop}
B_{k}(\partial M):=\{\sum_{i=1}^{k}\alpha_i\d_{a_i}, a_i\in \partial M, \alpha_i\geq 0, i=1,\cdots, k,\;\,\sum_{i=1}^{k}\alpha_i=k\},
\end{equation}
Finally, we set 
\begin{equation}\label{eq:defbarym}
A_{k, \bar k}:=\widetilde{B_{k}(\partial M)\times \bar B^{\bar k}_1},
\end{equation}
and 
\begin{equation}\label{eq:defbarynega}
A_{k-1, \bar k}:=\widetilde{B_{k-1}(\partial M)\times \bar B^{\bar k}_1}.
\end{equation}


\section{Blow-up analysis and critical points at infinity}\label{eq:cpi}
This section deals with the blowup analysis of sequences of vanishing viscosity solutions of the type  \begin{equation}\label{eq:blowupeq}
\left\{
\begin{split}
P^4_gu_l+2t_lQ_g&=0\;\;&\text{in}\;\;M,\\
P^3_gu_l+t_lT_g&=t_lKe^{3u}\;\;&\text{on}\;\;\partial M,\\
\frac{\partial u_l}{\partial n_g}&=0\;\;&\text{on}\;\;\partial M.
\end{split}
\right.
\end{equation}
with \;$t_l\rightarrow 1$\; under the assumption $\ker \mathbb{P}_g^{4, 3}\simeq \R$ and $\kappa_{(P^4, P^3)}=4k\pi^2$\; with \;$k\geq 1$ and their use to characterize the critical points at infinity of \;$\mathcal{E}_g$.
\subsection{Blow-up analysis}
The local behaviour of blowing up sequences of solutions of \eqref{eq:blowupeq} is understood. In fact, in \cite{nd6}, we prove the following lemma.
\begin{lem}\label{eq:local}
Assuming that $(u_{l})$ is a blowing up sequence of solutions to \eqref{eq:blowupeq}, then up to a subsequence, there exists \;$k$\;converging sequence of points\; $(x_{i,l})_{l\in \N}, x_{i, l}\in \partial M$\;\;with limits \;$x_i\in \partial M$, $i=1,\cdots,k $,\;\;$k$\; sequences $(\mu_{i,l})_{l\in \N} \;\;i=1,\cdots, k$\;of positive real numbers converging to \;$0$\;such that the following hold: \\\\
a)
\\
$$
\hspace{-45pt}\frac{d_{\hat g}(x_{i,l},x_{j,l})}{\mu_{i,l}}\longrightarrow +\infty \;\;\;\; i\neq j \;\;i,j =1, \cdots, k\;\;\;\;\;and \;\;\;\; t_{l}K(x_{i,l})\mu_{i.l}^{3}e^{3u_{l}(x_{i,l})}e^{-3\log 2}=2.
$$
b)
\\
$$
\hspace{-45pt}v_{i,l}(x)=u_{l}(exp^g_{x_{i,l}}(\mu_{i,l}x))-u_{l}(x_{i,l})+\log 2\longrightarrow V_{0}(x)\;\;\;\; in\; \;\;C^{4}_{loc}(\R^4_+),
$$
$$
 V_{0}|_{\R^3}(x):=\log\left(\frac{2}{1+|x|^{2}}\right).
$$
c)
\\
$$\hspace{-45pt}\text{There exists}\;\;C>0\;\; \text{such that}\;\;
\inf_{i=1,\cdots, k}d_{\hat g}(x_{i,l},x)^{3}e^{3u_{l}(x)}\leq C \;\;\;\;\forall x\in \partial M,\;\;\forall l\in \N.
$$
d)\\
\begin{equation*}
\begin{split}
t_lKe^{3u_l}dS_g\rightarrow 4\pi^2\sum_{i=1}^k\d_{x_{i}}  \;\;\;\text{in the sense of measure},\;\;\;\text{and}\\\lim_{l\rightarrow +\infty}\oint_{\partial M}t_lKe^{3u_l}dS_g=4\pi^2k.
\end{split}
\end{equation*}
e)\\
$$
u_l-\ov {(u_l)}_{Q, T}\rightarrow \sum_{i=1}^k G(x_{i}, \cdot)\;\;\;\text{in}\;\;\;C^3_{loc}(\partial M - \{x_1,\cdots, x_k\}),\;\;\;\;\ov{ (u_l)}_{Q, T}\rightarrow -\infty.
$$
\end{lem}
\vspace{8pt}

\noindent
As a Liouville type problem, the following Harnack type inequality is sufficient to get the global description of blowing up sequences of solutions needed to describe the critical points at infinity of \;$\mathcal{E}_g$.
\begin{pro}\label{eq:refinedestimate}
Assuming that \;$u_l$\; is a blowing up sequence of solutions to \eqref{eq:blowupeq}, then Lemma \ref{eq:local} holds, and keeping the notations in Lemma \ref{eq:local}, we have that the points \;$x_{i, l}$\; are uniformly isolated, namely there exists $0<\eta_k<\frac{\varrho_0}{10}$ (where $\varrho_0$ is as in \eqref{eq:varro0}) such that for \;$l$\; large eneough, there holds
\begin{equation}\label{eq:unifiso}
d_{\hat g}(x_{i, l}, x_{j, l})\geq 4\ov C\eta_k ,\;\;\;\forall i\neq j=1, \cdots, k.
\end{equation}
Moreover,  the scaling parameters \;$\l_{i, l}:=\mu_{i, l}^{-1}$\; are comparable, namely there exists a large positive constant \;$\Lambda_0$\; such that
\begin{equation}\label{eq:uniformiso}
\Lambda_0^{-1}\l_{j,l}\leq\l_{i,l} \leq \Lambda_0\l_{j,l},\;\;\forall\;i,j
\end{equation}
Furthermore, we have that the following estimate around the blow up points holds
\begin{equation}\label{eq:sharpestimate}
u_l(y)+\frac{1}{3}\log\frac{t_t K_l(x_{i, l})}{2}=\log\frac{2\l_{i, l}}{1+\l_{i, l}^2(d_{\hat g_{x_{i, l}}}(y, x_{i, l}))^2}+ O(d_{\hat g}(y, x_{i, l})),\;\;\;\;\forall\;y\in\;B^{\hat g}_{x_{i, l}}(\eta).
\end{equation}
\end{pro}
\vspace{6pt}

\noindent
To prove Proposition \ref{eq:refinedestimate}, as it is standard for Louiville type problems, one starts with the uniform isolation of blowing-up points. Indeed, we have 
\begin{lem}\label{eq:isolated}
Assuming that \;$(u_l)_{l\in \N}$ is a bubbling sequence of solutions to BVP \eqref{eq:blowupeq}, then keeping the notations in Lemma \ref{eq:local}, we have that the points \;$x_{i, l}$\; are uniformly isolated, namely there exists $0<\eta_k<\frac{\varrho_0}{10}$ (where $\varrho_0$ is as in \eqref{eq:varro0}) such that for \;$l$\; large eneough, there holds
\begin{equation}\label{eq:unifiso}
d_{\hat g}(x_{i, l}, x_{j, l})\geq 4\ov C\eta_k ,\;\;\;\forall i\neq j=1, \cdots, k.
\end{equation}
\end{lem}
\begin{pf}
The proof use the  integral method of Step 4 in \cite{nd1} and hence we will be skectchy in many arguments. As in \cite{nd1}, we first fix $\frac{1}{3}<\nu<\frac{2}{3}$, and for $i=1, \cdots, k$, we set
$$
\bar u_{i, l}(r)=Vol_{\hat g}(\partial B_{x_{i}}(r))^{-1}\int_{\partial B_{x_{i}}(r)}u_l(x)d\sigma_{\hat g}(x),\;\;\;\forall\;0\leq r<inj_{\hat g}(\partial M),
$$
and
$$
\psi_{i,l}(r)=r^{4\nu}exp(4\bar u_{i,l}(r)),\;\;\;\forall\;0\leq r<inj_g(\partial M).
$$
Furthermore, as in \cite{nd1}, we define \;$r_{i, l}$ as follows
\begin{equation}\label{eq:ril}
r_{i, l}:=\sup\{R_{\nu}\mu_{i, l}\leq r\leq \frac{R_{i, l}}{2}\;\;\text{such that}\;\;\psi_{i,l}^{'}(r)<0\;\;\text{in}\;[R_{\nu}\mu_{i, l}, r[\};
\end{equation}
where \;$R_{i, l}:=\min_{j\neq i}d_{\hat g}(x_{i, l}, x_{j, l})$. 
Thus, by continuity and the definition of \;$r_{i, l}$, we have that
\begin{equation}\label{eq:varr}
\psi_{i, l}^{'}(r_{i, l})=0
\end{equation}
Now, as in \cite{nd1}, to prove \eqref{eq:unifiso}, it suffices to show that \;$r_{i, l}$\; is bounded below by a positive constant in dependent of \;$l$. Thus, we assume by contradiction that (up to a subsequence) \;$r_{i, l}\rightarrow 0$\; as \;$l\rightarrow +\infty$\; and look for a contradiction. In order to do that, we use the integral representation formula \eqref{eq:G4integral} and argue as in Step 4 of \cite{nd1} to derive the following estimate
$$
\psi^{'}_{i,l}(r_{i,l})\leq (r_{i,l})^{3\nu-1}exp(\bar u_{i,l}(r_{i,l}))\left (3\nu-2C+o_{l}(1)+O_l(r_{i,l})\right).
$$
with \;$C>1$. So from\;$\frac{1}{3}<\nu<\frac{2}{3}$, $C>1$ and \;$r_{i, l}\longrightarrow 0$\; as \;$l\rightarrow +\infty$, we  deduce that for \;$l$\; large enough, there holds
\begin{equation}\label{eq:negativeb}
\psi_{i,l}^{'}(r_{i,l})<0.
\end{equation}
Thus, \eqref{eq:varr} and \eqref{eq:negativeb} lead to a contradiction, thereby concluding the proof of \eqref{eq:unifiso}. Hence, the proof of the Lemma is complete.
\end{pf}
\vspace{4pt}

\noindent
The next step to derive Proposition \ref{eq:refinedestimate} is to establish its weak $O(1)$ -version.
\begin{lem}\label{eq:roughestimate}
Assuming that $(u_l)_{l\in \N}$ is a bubbling sequence of solutions to BVP \eqref{eq:blowupeq}, then keeping the notations in Lemma \ref{eq:local} and Lemma \ref{eq:isolated}, we have that for \;$l$ large enough, there holds 
\begin{equation}\label{eq:classupinf}
u_l(x)+\frac{1}{3}\log\frac{t_l K(x_{i})}{2}=\log\frac{2\l_{i, l}}{1+\l_{i, l}^2(d_{\hat g_{x_i}}(x, x_{i}))^2}+ O(1),\;\;\;\;\forall\;x\in\;B^{x_i}_{x_{i}}(\eta_k),
\end{equation}
up to choosing \;$\eta_k$ smaller than in Lemma \ref{eq:isolated}.
\end{lem}
\begin{rem}\label{eq:comparability}
We point out that the comparability of the scaling parameters \;$\l_{i,l}$'s follows directly from Lemma \ref{eq:roughestimate}.
\end{rem}
\begin{pf}
We are going to use the method of \cite{nd7}, hence we will be sketchy in many arguments. Like in \cite{nd7}, thanks to Lemma \ref{eq:isolated}, we  will focus only on one blow-up point and called it $x\in \partial M$. Thus, we are in the situation where there exists  a sequence $x_l\in \partial M$\; such that \;$x_l\rightarrow x$\; with \;$x_l$\; local maximum point for \;$u_l$\; on $\partial M$\; and \;$u_l(x_l)\rightarrow +\infty$.
Now, we recall $g_x=e^{2u_x}g$ and choose $\eta_1$ such that
\;$20\eta_1<\min\{inj_g(\partial M), \varrho_0, \varrho_{k},  d\}$ with $4d\leq r_{i,l}$\; where \;$r_{i,l}$ is as in the proof of Lemma \ref{eq:isolated}. Next, we let \;$\hat w_x$ \;be the unique solution of the following boundary value problem
\begin{equation}\label{eq:bvpauxii}
\left\{
\begin{split}
P^4_{g_x}\hat w_x&=P^4_{\hat g}u_x\;\;&\text{in}\;\;M,\\
P^3_{g_x}\hat w_x&=P^3_{\hat g}u_x\;\;&\text{on}\;\;\partial M,\\
\frac{\partial \hat w_x}{\partial n_{g_x}}&=0\;\;&\text{on}\;\;\partial M,\\
\ov{\hat w}_{(Q, T)}&=0.
\end{split}
\right.
\end{equation}
Using standard elliptic regularity theory and \eqref{eq:proua}, we derive
\begin{equation}\label{eq:noinflu}
\hat w(y)=O(d_g(y, x)) \;\;\;\text{in}\;\;\;\; B_{x}^{g_x, +}(2\eta_1).
\end{equation}
On the other hand, using the conformal covariance properties of the Paneitz operator and of the Chang-Qing one, see \eqref{eq:law}, we have that \;$\hat u_l:=u_l-\hat w_x$\; satisfies
\begin{equation*}
\left\{
\begin{split}
P^4_{g_x}\hat u_l+2\hat Q_l&=0\;\;&\text{in}\;\;M,\\
P^3_{g_x}\hat u_l+\hat T_l&=t_lKe^{3\hat u_l}\;\;&\text{on}\;\;\partial M,\\
\frac{\partial \hat u_l}{\partial n_{g_x}}&=0\;\;&\text{on}\;\;\partial M.
\end{split}
\right.
\end{equation*}
with
$$
\hat Q_l=t_le^{-4\hat w}Q_g+\frac{1}{2}P^4_{\hat g}\hat w\;\;\;\text{and}\;\;\;\hat T_l=t_le^{-3\hat w}T_g+P^3_{\hat g}\hat w.
$$
Next, as in \cite{nd7}, we are going to establish the classical sup+inf-estimate for $\hat u_l$ , since thanks \eqref{eq:noinflu} all terms coming from \;$\hat w_x$\; can be absorbed on the right hand side of \eqref{eq:classupinf}.
Now, we are going to rescale the functions \;$\hat u_l$\; around the points \;$x$. In order to do that, we define \; $\varphi_l: B^{\R^3}_0(2\eta_1\mu_l^{-1})\longrightarrow B^{\hat g_x}_{x}(2\eta_1)$\;  by the formula $\varphi_l(z):=\mu_lz$\; and \;$\mu_l$\; is the corresponding scaling parameter given by Lemma \ref{eq:local}. Furthermore, as in \cite{nd7}, we define the following rescaling of \;$\hat u_l$ \;
$$
v_l:=\hat {u_l}\circ\varphi_l+\log \mu_l+\frac{1}{3}\log\frac{t_l K(x)}{2}.
$$
Using the Green's representation formula for and the method of the method of \cite{nd7}, we get
\begin{equation}\label{eq:estint}
v_l(z)+2\log|z|=O(1), \;\;\text{for}\;\;z\in \bar B^{\R^3}_0(\frac{\eta_1}{\mu_l})-B^{\R^3}_0(-\log \mu_l).
\end{equation}
Now, we are going to show that the estimate \eqref{eq:estint} holds also in \;$\bar B^{\R^3}_0(-\log \mu_l)$. To do so, we use Lemma \ref{eq:local}and the same arguments as in \cite{nd7} to deduce
\begin{equation}\label{eq:euclideanversion}
v_l(z)+2\log|z|=O(1), \;\;\text{for}\;\;\;z\in \bar B^{\R^3}_0(-\log \mu_l).
\end{equation}
Now, combining \eqref{eq:estint} and \eqref{eq:euclideanversion}, we obtain
\begin{equation}\label{eq:eucball}
v_l(z)+2\log|z|=O(1), \;\;\text{for}\;\;\;z\in \bar B^{\R^3}_0(\frac{\eta_1}{\mu_l}).
\end{equation}
Thus scaling back, namely using \;$y=\mu_lz$\; and the definition of \;$v_l$, we obtain the desired \;$O(1)$-estimate.  Hence the proof of the Lemma is complete.
\end{pf}

\begin{pfn}{ of formula \eqref{eq:sharpestimate} of Proposition \ref{eq:refinedestimate}}\\
We are going to use the method of \cite{nd7}, hence we will be sketchy in many arguments. Now, let $V_0$ be the unique solution of the following conformally invariant integral equation
$$
V_0(z)=\frac{1}{2\pi^2}\int_{\R^3}\log\frac{|y|}{|z-y|}e^{3V_0(y)}dy+\log 2,\;\;\; V_0(0)=\log 2, \n V_0(0)=0.
$$
Next, we set $w_l(z)=v_l(z)-V_0(z)$ for $z\in B^{\R^3}_0(\eta_1\mu_l^{-1})$, and use Lemma \ref{eq:roughestimate} \;to infer that
\begin{equation}\label{eq:weakes}
|w_l|\leq C\;\;\;\,\text{in}\;\;\;\; B^{\R^3}_0(\eta_1\mu_l^{-1})
\end{equation}
On the other hand, it is easy to see that to achieve our goal, it is sufficient to show
\begin{equation}\label{eq:3.2}
|w_l|\leq C\mu_l|z|\,\,\,\,\text{in}\,\,\,\, B^{\R^3}_0(\eta_1\mu_l^{-1}).
\end{equation}
To show \eqref{eq:3.2}, we first set
$$
\L_l:=\max_{z\in \Omega_l}\frac{|w_l(z)|}{\mu_l(1+|z|)}
$$
with
$$
\Omega_l=\ov B^{\R^3}_0(\eta_1\mu_l^{-1})
$$
We remark that to show \eqref{eq:3.2}, it is equivalent to prove that $\L_l$ is bounded.  Now, let us suppose that $\L_l\rightarrow +\infty$ as $l\rightarrow +\infty$, and look for a contradiction. To do so, we will use the method of  \cite{nd7}. For this, we first choose  a sequence of points $z_l\in \Omega_l$ such that $\L_l=\frac{|w_l(z_l)|}{\mu_l(1+|z_l|)}$. Next, up to a subsequence, we have that either $z_l\rightarrow z^{*}$ as $l\rightarrow+\infty$ (with $z^{*}\in \R^3$) or $|z_l|\rightarrow +\infty$ as $l\rightarrow+\infty$. Now, we make the following definition
$$
\bar w_l(z):=\frac{w_l(z)}{\L_l\mu_l(1+|z_l|)},
$$
and have
\begin{equation}\label{eq:3.3}
|\bar w_l(z)|\leq \left(\frac{1+|z|}{1+|z_l|}\right),
\end{equation}
and
\begin{equation}\label{eq:contradictionc}
|\bar w_l(z_l)|=1.
\end{equation}
Now, we consider the case where the points \;$z_l$\; escape to infinity.\\\\
{\em Case 1 }: $|z_l|\rightarrow +\infty$\\
In this case, using the integral representation \eqref{eq:G4integral} with respect to  \;$g_x$\; and the method of \cite{nd7}, we obtain
$$
\bar w_l(z_l)=\frac{1}{2\pi^2}\int_{\Omega_l}\log \frac{|\xi|}{|z_l-\xi|}\left(\frac{O(1)(1+|\xi|)^{-5}}{(1+|z_l|)}+\frac{O(1)(1+|\xi|)^{-5}}{\L_l(1+|z_l|)}\right)d\xi+o(1).
$$
Now, using the fact that $|z_l|\rightarrow+\infty$ as $l\rightarrow +\infty$, one can easily check that
$$
\bar w_l(z_l)=\frac{1}{2\pi^2}\int_{\Omega_l}\log \frac{|\xi|}{|z_l-\xi|}\left(\frac{O(1)(1+|\xi|)^{-5}}{(1+|z_l|)}+\frac{O(1)(1+|\xi|)^{-5}}{\L_l(1+|z_l|)}\right)d\xi=o(1).
$$ 
Hence, we reach a contradiction to \eqref{eq:contradictionc}. \\
Now, we are going to show that, when the points $z_l\rightarrow z^{*}$ as $l\rightarrow +\infty$, we reach a contradiction as well.\\\\
{\em Case 2}: $z_l\rightarrow z^{*}$\\
In this case, using the assumption $z_l\longrightarrow z^*$, the Green's representation formula, and the method of \cite{nd7}, we obtain that up to a subsequence
\begin{equation}\label{eq:lisupinf}
\bar w_l\rightarrow w \;\;\;\text{in}\;\;C^1_{loc}(\R^3) \;\;\;\text{as}\;\;\; l\rightarrow +\infty,\;\;\:\
\end{equation}
and
\begin{equation}\label{eq:barwl1b}
\begin{split}
\bar w_l(z)=\frac{1}{2\pi^2}\int_{\Omega_l} \log\frac{|\xi|}{|z-\xi|}\frac{K\circ\varphi_l(\xi)}{K\circ\varphi_l(0)}e^{3\vartheta_l(\xi)}\bar w(\xi)d\xi+\frac{1}{\L_l\mu_l(1+|z_l|)2\pi^2}\int_{\Omega_l} \log\frac{|\xi|}{|z-\xi|}O(\mu_l(1+|\xi|)^{-5})d\xi\\+\frac{O(1)+O(|z|)}{\L_l(1+|z_l|)},
\end{split}
\end{equation}
where $e^{3\theta_l}:=\int_0^1e^{3(sv_l+(1-s)V_0)}ds$.
Thus, appealing to \eqref{eq:lisupinf} and \eqref{eq:barwl1b}, we infer that $w$ satisfies
\begin{equation}\label{eq:wlim1}
w(z)=\frac{1}{2\pi^2}\int_{\R^3}\log\frac{|\xi|}{|z-\xi|}e^{3V_0(\xi)}w(\xi)d\xi
\end{equation}
Now, using \eqref{eq:3.3}, we have that $w$ satisfies the following asymptotics
\begin{equation}\label{eq:continf}
|w(z)|\leq C(1+|z|).
\end{equation}
On the other hand, from the definition of $v_l$,  it is easy to see that
\begin{equation}\label{eq:vanord1}
w(0)=0, \;\;\;\text{and}\;\;\;\;\n w(0)=0.
\end{equation}
So, using \eqref{eq:wlim1}-\eqref{eq:vanord1}, and observing that Lemma 3.7 in \cite{nd7} holds for dimension \;$3$, , we obtain
$$
w=0.
$$
 However, from \eqref{eq:contradictionc}, we infer that $w$ satisfies also
 \begin{equation}\label{eq:nonzero}
 |w(z^*)|=1
 \end{equation}
So we reach a contradiction in the second case also. Hence the proof of the lemma is complete.
\end{pfn}
\vspace{6pt}

\noindent
Because of the lack of understanding of the blowing PS-sequences for Louiville type problems, the role of the PS-sequences can be replaced by the vanishing viscosity solutions of the type of \eqref{eq:blowupeq} via the  following Bahri-Lucia's deformation lemma.
\begin{lem}\label{eq:deformlem}
 Assuming that \;$a, \;b\in \R$\; such that \;$a<b$ and  there is no critical values of \;$\mathcal{E}_g$\; in \;$[a, b]$, then there are two possibilities\\
1) Either  $$(\mathcal{E}_g)^a\;\text{is a deformation retract of}\;\; (\mathcal{E}_g)^b.$$
2) Or there exists a sequence \;$t_l\rightarrow 1$\; as \;$l\rightarrow +\infty$\; and a sequence of critical point \;$u_l$\; of \;$(\mathcal{E}_g)_{t_l}$\; verifying \;$a\leq \mathcal{E}_g(u_l)\leq b$\; for all \;$l\in \N^*$, where \;$(\mathcal{E}_g)_{t_l}$\; is as in \eqref{eq:jt} with $t$ replaced by \;$t_l$.
\end{lem}
\vspace{6pt}

\noindent
On the other hand, setting
\begin{equation}\label{eq:ninfinity}
\begin{split}
V_R(k, \epsilon, \eta):=\{u\in \mathcal{H}_{\frac{\partial }{\partial n}}:\;\;\exists a_1, \cdots, a_k\in \partial M, \;\;\l_1,\cdots, \l_k>0, \;\;||u-\ov{u}_{Q, T}-\sum_{i=1}^k\varphi_{a_i,\lambda_i}||_{\mathbb{P}^{4, 3}}=\\O\left( \sum_{i=1}^k\frac{1}{\l_i}\right)\;\;\;\lambda_i\geq \frac{1}{\epsilon},\;\;\frac{2}{\Lambda}\leq \frac{\l_i}{\l_j}\leq \frac{\Lambda}{2},\;\;\text{and}\;\;d_{\hat g}(a_i, a_j)\geq 4\ov C\eta\;\;\text{for}\;i\neq j\},
\end{split}
\end{equation}
where \;$\ov C$\; is as in \eqref{eq:proua}, \;$L$\; as in \eqref{eq:lninfinity}, , $O(1):=O_{A, \bar \l, u, \epsilon}(1)$ meaning bounded uniformly in \;$\bar\l:=(\l_1, \cdots, \l_k)$, $A:=(a_1, \cdots, a_k)$, $u$, $\epsilon$, we have as in \cite{nd7} that Proposition \ref{eq:refinedestimate} implies the following one.
\begin{lem}\label{eq:escape}
Let $\epsilon$ and $\eta$ be small positive real numbers with \;$0<2\eta<\varrho$ where \;$\varrho$\; is as in \eqref{eq:cutoff}. Assuming that $u_l$ is a sequence of blowing up critical point  of \;$(\mathcal{E}_g)_{t_l}$\; with \;$\ov{(u_l)}_{Q, T}=0, l\in \N$\; and \;$t_l\rightarrow 1$\; as \;$l\rightarrow +\infty$ , then  there exists \;$l_{\epsilon, \eta}$\; a large positive integer such that for every \;$l\geq l_{\epsilon, \eta}$, we have \;$u_l\in V_R(k, \epsilon, \eta)$, and for the definition of \;$V_R(k, \epsilon, \eta)$, see \eqref{eq:ninfinity}.
\end{lem}
\vspace{6pt}

\noindent
Finally, as in \cite{nd7}, we have that Lemma \ref{eq:deformlem} and Lemma \ref{eq:escape} implies the following one.
\begin{lem}\label{eq:deformlemr}
Assuming that \;$\epsilon$\; and \;$\eta$\; are small positive real numbers with \;$0<2\eta<\varrho$, then for\; $a, b\in \R$ such that $a<b$, we have that if  there is no critical values of \;$\mathcal{E}_g$\; in \;$[a, b]$, then there are two possibilities\\
1) Either  $$(\mathcal{E}_g)^a\;\text{is a deformation retract of}\;\; (\mathcal{E}_g)^b.$$
2) Or there exists a sequence \;$t_l\rightarrow 1$\; as $l\rightarrow +\infty$\; and a sequence of critical point \;$u_l$\; of \;$(\mathcal{E}_g)_{t_l}$\; (for its definition see \eqref{eq:jt}) verifying \;$a\leq \mathcal{E}_g(u_l)\leq b$\; for all \;$l\in \N^*$\; and \;$l_{\epsilon, \eta}$\; a large positive integer such that  \;$u_l\in V_R(k, \epsilon, \eta)$\; for all \;$l\geq l_{\epsilon, \eta}$,  and for the definition of $V_R(k, \epsilon, \eta)$, see \eqref{eq:ninfinity}.
\end{lem}


\subsection{Energy and gradient estimates at infinity}
In this subsection, we present energy and gradient estimates needed to characterize the critical points at infinity of \;$\mathcal{E}_g$. We start with a parametrization of infinity. Indeed, as a Liouville type problem, we have that for \;$\eta$\; a small positive real number with \;$0<2\eta<\varrho$, there exists \;$\epsilon_0=\epsilon_0(\eta)>0$ such that $\forall\;0<\epsilon\leq \epsilon_0$, we have 
\begin{equation}\label{eq:mini}
\forall u\in V_R(k, \epsilon, \eta), \\\text{the minimization problem }\min_{B_{\epsilon, \eta}}||u-\ov{u}_{Q, T}-\sum_{i=1}^k\alpha_i\varphi_{a_i, \l_i}-\sum_{r=1}^{\bar k}\beta_r(v_r-\ov{(v_r)}_{Q, T})||_{\mathbb{P}^{4, 3}}
\end{equation}
has a unique solution, up to permutations, where \;$B_{\epsilon, \eta}$\; is defined as follows
\begin{equation}
\begin{split}
{B_{\epsilon, \eta}:=\{(\bar\alpha, A, \bar \l, \bar \beta)\in \R^k\times (\partial M)^k\times (0, +\infty)^k\times \R^{\bar k}}:|\alpha_i-1|\leq \epsilon, \l_i\geq \frac{1}{\epsilon}, i=1, \cdots, k, \\d_{\hat g}(a_i, a_j)\geq 4\ov C\eta, i\neq j, |\beta_r|\leq R, r=1, \cdots, \bar k\}.
\end{split}
\end{equation}
 Moreover, using the solution of \eqref{eq:mini}, we have that every \;$u\in V_R(k, \epsilon, \eta)$\; can be written as
\begin{equation}\label{eq:para}
u-\ov{u}_{Q, T}=\sum_{i=1}^k\alpha_i\varphi_{a_i, \l_i}+\sum_{r=1}^{\bar k}\beta_r(v_r-\ov{(v_r)}_{Q, T})+w,
\end{equation}
where \;$w$\; verifies the following orthogonality conditions
\begin{equation}\label{eq:ortho}
\begin{split}
\ov {w}_{Q, T}=\left\langle\varphi_{a_i, \l_i}, w\right\rangle_{\mathbb{P}^{4, 3}}=\left\langle\frac{\partial\varphi_{a_i, \l_i}}{\partial \l_i}, w\right\rangle_{\mathbb{P}^{4, 3}}=\left\langle\frac{\partial\varphi_{a_i, \l_i}}{\partial a_i}, w\right\rangle_{\mathbb{P}^{4, 3}}=\left\langle v_r, w\right\rangle_{\mathbb{P}^{4, 3}}=0, i=1, \cdots, k,\\r=1, \cdots, \bar k
\end{split}
\end{equation}
and the estimate
\begin{equation}\label{eq:estwmin}
||w||_{\mathbb{P}^{4, 3}}=O\left(\sum_{i=1}^k\frac{1}{\l_i}\right),
\end{equation}
where here $O\left(1\right):=O_{\bar \alpha, A, \bar \l, \bar\beta , w, \epsilon}\left(1\right)$. Furthermore, the concentration points \;$a_i$,  the masses \;$\alpha_i$, the concentrating parameters \;$\l_i$\; and the negativity parameter \;$\beta_r$ \;in \eqref{eq:para} verify also
\begin{equation}\label{eq:afpara}
\begin{split}
d_{\hat g}(a_i, a_j)\geq 4\ov C\eta,\;i\neq j=1, \cdots, k, \frac{1}{\L}\leq\frac{\l_i}{\l_j}\leq\L\;\;i, j=1, \cdots, k, \;\;\l_i\geq\frac{1}{\epsilon},\;\;\text{and}\\\;\;\;\sum_{r=1}^{\bar k}|\beta_r|+\sum_{i=1}^k|\alpha_i-1|\sqrt{\log \l_i}=O\left(\sum_{i=1}^k\frac{1}{\l_i}\right)
\end{split}
\end{equation}
with still $O\left(1\right)$ as in \eqref{eq:estwmin}.
\vspace{6pt}

\noindent
Because of the translation invariant property of \;$\mathcal{E}_g$\; and the parametrization \eqref{eq:para}, to derive energy estimate in \;$V_R(k, \epsilon, \eta)$\; we start with the following lemma.
\begin{lem}\label{eq:energyest}
Assuming that \;$\eta$\; is a small positive real number with \;$0<2\eta<\varrho$\; where \;$\varrho$\; is as in \eqref{eq:cutoff}, and \;$0<\epsilon\leq \epsilon_0$\; where \;$\epsilon_0$\; is as in \eqref{eq:mini}, then  for \;$a_i\in M$\; concentration points,  $\alpha_i$\; masses, $\l_i$\; concentration parameters ($i=1,\cdots,k$), and \;$\beta_r$ negativity parameters ($r=1, \cdots, \bar k$) satisfying \eqref{eq:afpara}, we have
\begin{equation*}
\begin{split}
&\mathcal{E}_g\left(\sum_{i=1}^k\alpha_i\varphi_{a_i, \l_i}+\sum_{r=^1}^{\bar k}\beta_r(v_r-\ov{(v_r)}_{Q, T})\right)=C^{k}_0-8\pi^2\mathcal{F}_K(a_1, \dots, a_k)+2\sum_{r=1}^{\bar k}\mu_r\beta_r^2\\&+\sum_{i=1}^k(\alpha_i-1)^2\left[16\pi^2\log\l_i+8\pi^2H(a_i, a_i)+C^{k}_1\right]\\&+8\pi^2\sum_{i=1}^k(\alpha_i-1)\left[\sum_{r=1}^{\bar k}2\beta_r(v_r-\ov{(v_r)}_{Q, T})(a_i)+\sum_{ j=1, j\neq i}^k (\alpha_j-1)G(a_i, a_j)\right]\\&-\frac{c^18\pi^2}{9}\sum_{i=1}^k\frac{1}{\l_i^2}\left(\frac{\D_{\hat g_{a_i}}\mathcal{F}^{A}_i(a_i)}{\mathcal{F}^{A}_i(a_i)}-\frac{3}{4}R_{\hat g}(a_i)\right)\\&+\frac{c^18\pi^2}{9}\sum_{i=1}^k\frac{\tilde \tau_i}{\l_i^2}\left(\frac{\D_{\hat g_{a_i}}\mathcal{F}^{A}_i(a_i)}{\mathcal{F}^{A}_i(a_i)}-\frac{3}{4}R_{\hat g}(a_i)\right)\\&+\frac{16\pi^2}{3}\sum_{i=1}^k\log(1-\tilde\tau_i)+O\left(\sum_{i=1}^k|\alpha_i-1|^3+\sum_{r=1}^{\bar k}|\beta_r|^3+\sum_{i=1}^k\frac{1}{\l^3_i}\right),
\end{split}
\end{equation*}
where \;$O\left(1\right)$ means here \;$O_{\bar\alpha, A, \bar \l, \bar \beta, \epsilon}\left(1\right)$ \;with \;$\bar\alpha=(\alpha_1, \cdots, \alpha_k)$, $A:=(a_1, \cdots, a_k)$, \;$\bar \l:=(\l_1, \cdots, \l_k)$, $\bar \beta:=(\beta_1, \cdots, \beta_{\bar k})$\; and for \;$i=1, \cdots, k$, $$\tilde \tau_i:=1-\frac{k\tilde\gamma_i}{\Gamma},\;\;\;\; \Gamma:=\sum_{i=1}^k\tilde\gamma_i, \;\;\;\tilde\gamma_i:=\tilde c_i\l_i^{6\alpha_i-3}\mathcal{F}^{A}_i(a_i)\mathcal{G}_i(a_i),$$ with 
$$
\tilde c_i:=\int_{\R^3}\frac{1}{(1+|y|^2)^{3\alpha_i}}dy
$$
\begin{equation*}
 \begin{split}
\mathcal{G}_i(a_i):=e^{3((\alpha_i-1)H(a_i, a_i)+\sum_{j=1, j\neq i}^k(\alpha_j-1)G(a_j, a_i))}e^{\frac{3}{2}\sum_{j=1, j\neq i}^k\frac{\alpha_j}{\l_j^2}\D_{g_{a_j}}G(a_j, a_i)}e^{\frac{3}{2}\frac{\alpha_i}{\l_i^2}\D_{g_{a_i}}H(a_i, a_i)}\\\times e^{3\sum_{r=1}^{\bar k}\beta_rv_r(a_i)},
\end{split}
\end{equation*}
 $C_0^{k}$\; is a real number depending only on \;$k$, $C_1^k$\; is a real number \; and \;$c^1$\; is a positive real number and for the meaning of \;$O_{\bar\alpha, A, \bar \l, \bar \beta, \epsilon}\left(1\right)$.
\end{lem}
\begin{pf}
The proof is the same as the one Lemma 4.1 in \cite{no1} replacing Lemma 10.1- Lemma 10.4 in \cite{no1} by Lemma \ref{eq:intbubbleest}-Lemma \ref{eq:interactest}.
\end{pf}
\vspace{6pt}

\noindent
Concerning the gradient estimates  of \;$\mathcal{E}_g$\;  in \;$V_R(k, \epsilon, \eta)$, we have in the directions of the scaling parameters :
\begin{lem}\label{eq:gradientlambdaest}
Assuming that \;$\eta$\; is a small positive real number with \;$0<2\eta<\varrho$\; where \;$\varrho$\; is as in \eqref{eq:cutoff}, and \;$\epsilon\leq \epsilon_0$\; where \;$\epsilon_0$\; is as in \eqref{eq:mini}, then  for \;$a_i\in \partial M$ concentration points,  $\alpha_i$\; masses , $\l_i$\; concentration parameters ($i=1,\cdots, k$) and \;$\beta_r$ negativity parameters ($r=1, \cdots, \bar k$) satisfying \eqref{eq:afpara}, we have that for every $r=1,\cdots, k$, there holds
\begin{equation*}
\begin{split}
&\left\langle\n^{\mathbb{P}^{4, 3}} \mathcal{E}_g\left(\sum_{i=1}^k\alpha_i\varphi_{a_i, \l_i}+\sum_{r=1}^{\bar k}\beta_r(v_r-\ov{(v_r)}_{Q, T})\right), \l_j\frac{\partial \varphi_{a_j, \l_j}}{\partial \l_j}\right\rangle_{\mathbb{P}^{4, 3}}= 16\pi^2\alpha_j\tau_j\\&-\frac{c^28\pi^2}{3\l_j^2}\left(\frac{\D_{\hat g_{a_j}} \mathcal{F}^{A}_j(a_j)}{\mathcal{F}^{A}_j(a_j)}-\frac{3}{4}R_{\hat g}(a_j)\right)-\frac{16\pi^2}{\l_j^2}\tau_j \D_{\hat g_{a_j}}H(a_j, a_j)\\&-\frac{16\pi^2}{\l_j^2}\sum_{i=1, i\neq j}^k\tau_i \D_{\hat g_{a_j}}G(a_j, a_i)+ \frac{c^28\pi^2}{\l_j^2}\tau_j\left(\frac{\D_{\hat g_{a_j}} \mathcal{F}^{A}_j(a_j)}{\mathcal{F}^{A}_j(a_j)}-\frac{3}{4}R_{\hat g}(a_j)\right)\\&+O\left(\sum_{i=1}^k|\alpha_i-1|^2+\sum_{r=1}^{\bar k}|\beta_r|^3+\sum_{i=1}^k\frac{1}{\l_i^3}\right),
\end{split}
\end{equation*}
where \;$A:=(a_1, \cdots, a_k)$, $O\left(1\right)$\; is as in Lemma \ref{eq:energyest}, $c^2$\; is a positive real number, and for \;$i=1, \cdots, k$, $$\tau_i:=1-\frac{k\tilde\gamma_i}{D}, \;\;\;\;\;D:=\oint_{\partial M} K(x)e^{3(\sum_{i=1}^k\alpha_i\varphi_{a_i, \l_i}(x)+\sum_{r=1}^{\bar k}\beta_r v_r(x))}dS_g(x),$$ with \;$\tilde\gamma_i$\; as in Lemma \ref{eq:energyest}.
\end{lem}
\begin{pf}
The proof is the same as the one Lemma 5.1 in \cite{no1} replacing Lemma 10.1- Lemma 10.4 in \cite{no1} by Lemma \ref{eq:intbubbleest}-Lemma \ref{eq:interactest}.
 \end{pf}
\vspace{6pt}
 
\noindent
 As in \cite{no1}, Lemma \ref{eq:gradientlambdaest} implies the following corollary.
 \begin{cor}\label{eq:cgradientlambdaest}
Assuming that \;$\eta$\; is a small positive real number with \;$0<2\eta<\varrho$\; where \;$\varrho$\; is as in \eqref{eq:cutoff},  and \;$0<\epsilon\leq \epsilon_0$\; where \;$\epsilon_0$\; is as in \eqref{eq:mini}, then for \;$a_i\in \partial M$\; concentration points,  $\alpha_i$\; masses, $\l_i$\; concentration parameters ($i=1,\cdots,k$), and \;$\beta_r$\; negativity parameters  ($r=1, \cdots, \bar k$)  satisfying \eqref{eq:afpara}, we have
\begin{equation*}
\begin{split}
&\left\langle\n^{\mathbb{P}^{4, 3}}\mathcal{E}_g\left(\sum_{i=1}^k\alpha_i\varphi_{a_i, \l_i}+\sum_{r=1}^{\bar k}\beta_r(v_r-\ov{(v_r)}_{(Q, T)})\right), \sum_{i=1}^k\frac{\l_i}{\alpha_i}\frac{\partial \varphi_{a_i, \l_i}}{\partial \l_i}\right\rangle_{\mathbb{P}^{4, 3}}=\\& \sum_{i=1}^k\frac{c^38\pi^2}{\l_i^2}\left(\frac{\D_{\hat g_{a_i}} \mathcal{F}^{A}_i(a_i)}{\mathcal{F}^{A}_i(a_i)}-\frac{3}{4}R_{\hat g}(a_i)\right)\\&+O\left(\sum_{i=1}^k|\alpha_i-1|^2+\sum_{r=1}^{\bar k}|\beta_r|^3+\sum^{k}_{i=1}\tau^3_i+\sum_{i=1}^k\frac{1}{\l_i^3}\right),
\end{split}
\end{equation*}
where \;$A:=(a_1, \cdots, a_k)$, $O\left(1\right)$ as as in Lemma \ref{eq:energyest}, $c^3$\; is a positive real number, and for \;$i=1, \cdots, k$, \;$\tau_i$\; is as in Lemma \ref{eq:gradientlambdaest}.
\end{cor}
\begin{pf}
The proof uses the strategy of the proof of Corollary 5.2 in \cite{no1} replacing Lemma 5.1 in \cite{no1}  by its counterpart Lemma \ref{eq:gradientlambdaest}.
\end{pf}
\vspace{6pt}

\noindent
For the gradient estimate in the directions of mass concentrations, we have:
\begin{lem}\label{eq:gradientalpha}
Assuming that \;$\eta$\; is a small positive real number with \;$0<2\eta<\varrho$\; where \;$\varrho$\; is as in \eqref{eq:cutoff}, and \;$0<\epsilon\leq \epsilon_0$\; where \;\;$\epsilon_0$\; is as in \eqref{eq:mini}, then  for \;$a_i\in \partial M$\; concentration points,  $\alpha_i$\; masses, $\l_i$\; concentration parameters ($i=1,\cdots, k$), and \;$\beta_r$\; negativity parameters ($r=1, \cdots, \bar k$) satisfying \eqref{eq:afpara}, we have that for every \;$j=1, \cdots, k$, there holds
\begin{equation*}
\begin{split}
&\left\langle\n^{\mathbb{P}^{4, 3}}\mathcal{E}_g\left(\sum_{i=1}^k\alpha_i\varphi_{a_i, \l_i}+\sum_{r=1}^{\bar k}\beta_r(v_r-\ov{(v_r)}_{(Q, T)})\right), \varphi_{a_j, \l_j}\right\rangle_{\mathbb{P}^{4, 3}}=\\&\left(2\log\l_j+H(a_j, a_j)-C_2\right)\frac{1}{\alpha_j}\left\langle\n^{\mathbb{P}^{4, 3}} \mathcal{E}_g\left(\sum_{i=1}^k\alpha_i\varphi_{a_i, \l_i}+\sum_{r=1}^{\bar k}\beta_r(v_r-\ov{(v_r)}_{(Q, T)})\right), \l_j\frac{\partial \varphi_{a_j, \l_j}}{\partial \l_j}\right\rangle_{\mathbb{P}^{4, 3}}\\&+\sum_{i=1, i\neq j}^k G(a_j, a_i)\left\langle\n^{\mathbb{P}^{4, 3}} \mathcal{E}_g\left(\sum_{i=1}^k\alpha_i\varphi_{a_i, \l_i}+\sum_{r=1}^{\bar k}\beta_r(v_r-\ov{(v_r)}_{(Q, T)})\right), \l_i\frac{\partial \varphi_{a_i, \l_i}}{\partial \l_i}\right\rangle_{\mathbb{P}^{4, 3}}\\&+32\pi^2(\alpha_j-1)\log\l_j+O\left(\log\l_j\left(\sum_{i=1}^k\frac{|\alpha_i-1|}{\log \l_i}+(\sum_{r=1}^{\bar k}|\beta_r|)(\sum_{i=1}^k\frac{1}{\log \l_i})+\sum_{i=1}^k\frac{1}{\l^2_i}\right)\right),
\end{split}
\end{equation*}
where \;$O\left(1\right)$\; as as in Lemma \ref{eq:energyest} and \;$C_2$\; is a real number.
\end{lem}
\begin{pf}
It follows from the same arguments as in Lemma 5.3 in \cite{no1}.
\end{pf}
\vspace{6pt}

\noindent
Concerning the gradient estimate in the directions of points of concentrations, we have:
\begin{lem}\label{eq:gradientaest}
Assuming that \;$\eta$\; is a small positive real number with \;$0<2\eta<\varrho$\; where \;$\varrho$\; is as in \eqref{eq:cutoff}, and \;$0<\epsilon\leq \epsilon_0$\; where $\epsilon_0$ is as in \eqref{eq:mini}, then  for \;$a_i\in \partial M$\; concentration points,  $\alpha_i$\; masses, $\l_i$\; concentration parameters ($i=1,\cdots, k$), and \;$\beta_r$\; negativity parameters ($r=1, \cdots, \bar k$) satisfying \eqref{eq:afpara}, we have that for every \;$j=1, \cdots, k$, there holds
\begin{equation*}
\begin{split}
\left\langle\n^{\mathbb{P}^{4, 3}} \mathcal{E}_g\left(\sum_{i=1}^k\alpha_i\varphi_{a_i, \l_i}+\sum_{r=1}^{\bar k}\beta_r(v_r-\ov{(v_r)}_{(Q, T)})\right), \frac{1}{\l_j}\frac{\partial \varphi_{a_j, \l_j}}{\partial a_j}\right\rangle_{\mathbb{P}^{4, 3}}&=-\frac{c^2 32\pi^2}{\l_j}\frac{\n_{\hat g}\mathcal{F}_j^{A}(a_j)}{\mathcal{F}_j^{A}(a_j)}\\&+O\left(\sum_{i=1}^k|\alpha_i-1|^2\right)\\&+O\left(\sum_{i=1}^k\frac{1}{\l_i^2}+\sum_{r=1}^{\bar k}|\beta_r|^2+\sum_{i=1}^k\tau_i^2\right),
\end{split}
\end{equation*}
where \;$A:=(a_1, \cdots, a_k)$, $O(1)$\; is as in Lemma \ref{eq:energyest}, $c^2$\; is as in Lemma \ref{eq:gradientlambdaest} and for  \;$i=1, \cdots, k$, \;$\tau_i$\; is as in Lemma \ref{eq:gradientlambdaest}.
\end{lem}
\begin{pf}
The proof is the same as the one of Lemma 5.4 in \cite{no1}.
\end{pf}
\vspace{6pt}

\noindent
Concerning the gradient estimate in the directions of the negativity parameters, we have:
\begin{lem}\label{eq:gradientbeta}
Assuming that $\eta$ is a small positive real number with \;$0<2\eta<\varrho$\; where \;$\varrho$\; is as in \eqref{eq:cutoff}, and \;$0<\epsilon\leq \epsilon_0$\; where \;$\epsilon_0$\; is as in \eqref{eq:mini}, then for \;$a_i\in \partial M$\; concentration points,  $\alpha_i$\; masses, $\l_i$ concentration parameters ($i=1,\cdots, k$), ad $\beta_r$ negativity parameters ($r=1, \cdots, \bar k$) satisfying \eqref{eq:afpara}, we have that for every \;$l=1, \cdots, \bar k$, there holds
\begin{equation*}
\begin{split}
\left\langle\n^{\mathbb{P}^{4, 3}} \mathcal{E}_g\left(\sum_{i=1}^k\alpha_i\varphi_{a_i, \l_i}+\sum_{r=1}^{\bar k}\beta_r (v_r-\ov{(v_r)}_{(Q, T)})\right), v_l-\ov{(v_l)}_{Q, T}\right\rangle_{\mathbb{P}^{4, 3}}=&4\mu_l\beta_l+O\left(\sum_{i=1}^k|\alpha_i-1|+\sum_{i=1}^k |\tau_i|\right)\\&+O\left(\sum_{i=1}^k\frac{1}{\l_i^2}\right),
\end{split}
\end{equation*}
where \;$O(1)\;$ is as in Lemma \ref{eq:energyest} and for  \;$i=1, \cdots, k$, \;$\tau_i$\; is as in Lemma \ref{eq:gradientlambdaest}
\end{lem}
\begin{pf}
It follows from the same arguments as in the proof of Lemma 5.5 in \cite{no1}.
\end{pf}

\subsection{Finite-dimensional reduction}
In this subsection, we complete the energy estimate of \;$\mathcal{E}_g$\; on \;$V_R(k, \epsilon, \eta)$\; via Lyapunov finite dimensional type reduction and second variation arguments. First of all, we have:
\begin{pro}\label{eq:expansionJ1}
Assuming that \;$\eta$\; is a small positive real number with \;$0<2\eta<\varrho$\; where \;$\varrho$\; is as in \eqref{eq:cutoff}, and \;$0<\epsilon\leq \epsilon_0$\; where \;$\epsilon_0$\; is as in \eqref{eq:mini} and \;$u=\ov{u}_{Q, T}+\sum_{i=1}^k\alpha_i\varphi_{a_i, \l_i}+\sum_{r=1}^{\bar k}\beta_r (v_r-\ov{(v_r)}_{(Q, T)})+w\in V_R(k, \epsilon, \eta)$ with \;$w$, the concentration points \;$a_i$,  the masses \;$\alpha_i$, the concentrating parameters \;$\l_i$  ($i=1, \cdots, k$), and the negativity parameters \;$\beta_r$ ($r=1, \cdots, \bar k$) verifying \eqref{eq:ortho}-\eqref{eq:afpara}, then we have
\begin{equation}\label{eq:exparoundbubble}
\mathcal{E}_g(u)=\mathcal{E}_g\left(\sum_{i=1}^k\alpha_i\varphi_{a_i, \l_i}+\sum_{r=1}^{\bar k}\beta_r (v_r-\ov{(v_r)}_{(Q, T)})\right)-f(w)+Q(w)+o(||w||_{\mathbb{P}^{4, 3}}^2),
\end{equation}
where
\begin{equation}\label{eq:linear}
f(w):=16\pi^2\frac{\oint_{\partial M} Ke^{3\sum_{i=1}^k\alpha_i\varphi_{a_i, \l_i}+3\sum_{r=1}^{\bar k}\beta_r v_r}wdS_g}{\oint_{\partial M} Ke^{3\sum_{i=1}^k\alpha_i\varphi_{a_i, \l_i}+3\sum_{r=1}^{\bar k}\beta_r v_r}dS_g},
\end{equation}
and
\begin{equation}\label{eq:quadratic}
Q(w):=||w||_{\mathbb{P}^{4, 3}}^2-24\pi^2\frac{\oint_{\partial M} Ke^{3\sum_{i=1}^k\alpha_i\varphi_{a_i, \l_i}+3\sum_{r=1}^{\bar k}\beta_r v_r}w^2dS_g}{\oint_{\partial M} Ke^{3\sum_{i=1}^k\alpha_i\varphi_{a_i, \l_i}+3\sum_{r=1}^{\bar k}\beta_r v_r}dS_g}.
\end{equation}
Moreover, setting
\begin{equation}\label{eq:eali}
\begin{split}
E_{a_i, \l_i}:=\{w\in \mathcal{H}_{\frac{\partial }{\partial n}}: \;\;\langle\varphi_{a_i, \l_i}, w\rangle_{\mathbb{P}^{4, 3}}=\langle\frac{\partial\varphi_{a_i, \l_i}}{\partial \l_i}, w\rangle_{\mathbb{P}^{4, 3}}=\langle\frac{\partial\varphi_{a_i, \l_i}}{\partial a_i}, w\rangle_{\mathbb{P}^{4, 3}}=0,\\\, \ov{w}_{(Q, T)}=\langle v_r, w\rangle_{\mathbb{P}^{4, 3}}=0, \;r=1, \cdots, \bar k,\,\;\text{and}\;\;||w||_{\mathbb{P}^{4, 3}}=O\left(\sum_{i=1}^k\frac{1}{\l_i}\right)\},
\end{split}
\end{equation}
and
\begin{equation}\label{eq:eal}
A:=(a_1, \cdots, a_k), \;\;\bar \l=(\l_1, \cdots, \l_k), \;\;E_{A, \bar \l}:=\cap_{i=1}^k E_{a_i, \l_i},
\end{equation}
we have that, the quadratic form \;$Q$\; is positive definite in \;$E_{A, \bar \l}$. Furthermore, the linear part \;$f$\; verifies that, for every \;$w\in E_{A, \bar \l}$, there holds
\begin{equation}\label{eq:estlinear}
f(w)=O\left[ ||w||_{\mathbb{P}^{4, 3}}\left(\sum_{i=1}^k\frac{|\n_{\hat g} \mathcal{F}^{A}_i(a_i)|}{\l_i}+\sum_{i=1}^k|\alpha_i-1|\log \l_i+\sum_{r=1}^{\bar k}|\beta_r|+\sum_{i=1}^k\frac{\log \l_i}{\l_i^{2}}\right)\right].
\end{equation}
where here \;$o(1)=o_{ \bar\alpha, A, \bar \beta, \bar\l, w, \epsilon}(1)$\; and \;$O\left(1\right):=O_{\bar\alpha, A, \bar \beta, \bar\l, w, \epsilon}\left(1\right)$.
\end{pro}
\vspace{6pt}

\noindent
As in \cite{no1}, to prove Proposition \ref{eq:expansionJ1}, we will need the following three coming lemmas. We start with the following one:
\begin{lem}\label{eq:holmos}
Assuming the assumptions of Proposition \ref{eq:expansionJ1} and \;$\gamma\in (0, 1)$\; small, then for every \;$q\geq 1$, there holds the following estimates\\
\begin{equation}\label{eq:esthm1}
\oint_{\partial M} Ke^{3\sum_{i=1}^k\alpha_i\varphi_{a_i, \l_i}+3\sum_{r=1}^{\bar k}\beta_r v_r}|w|^q=O\left(||w||_{\mathbb{P}^{4, 3}}^q(\sum_{i=1}^k\l_i^{3+\gamma})\right),
\end{equation}
\begin{equation}\label{eq:esthm2}
\oint_{\partial M} Ke^{3\sum_{i=1}^k\alpha_i\hat\d_{a_i, \l_i}+3\sum_{r=1}^{\bar k}\beta_r v_r}|w|^q=O\left(||w||_{\mathbb{P}^{4, 3}}^q(\sum_{i=1}^k\l_i^{\gamma})\right),
\end{equation}
\begin{equation}\label{eq:esthm3}
\oint_{\partial M} e^{3\hat\d_{a_i, \l_i}+3\sum_{r=1}^{\bar k}\beta_r v_r}d_{\hat g_{a_i}}(a_i, \cdot)|w|^q=O\left(||w||_{\mathbb{P}^{4, 3}}^q\frac{1}{\l_i^{1-\gamma}}\right),\;\; i=1, \cdots, k,
\end{equation}
\begin{equation}\label{eq:esthm4}
\oint_{\partial M} Ke^{3\sum_{i=1}^k\alpha_i \varphi_{a_i, \l_i}+3\sum_{r=1}^{\bar k}\beta_r v_r}e^{3\theta_w w}|w|^q=O\left(||w||_{\mathbb{P}^{4, 3}}^q(\sum_{i=1}^k\l_i^{3+\gamma})\right),
\end{equation}
where \;$\theta_w\in [0, 1]$, and
\begin{equation}\label{eq:esthm5}
\oint_{\partial M} Ke^{3\sum_{i=1}^k\alpha_i\varphi_{a_i, \l_i}+3\sum_{r=1}^{\bar k}\beta_r v_r}\left(e^{3w}-1-3w-\frac{9}{2}w^2\right)dV_g=o\left(||w||_{\mathbb{P}^{4, 3}}^2(\sum_{i=1}^k\l_i^{3})\right).
\end{equation}
where here \;$o(1)$\; and \;$O\left(1\right)$\; are as in Proposition \ref{eq:expansionJ1}.
\end{lem}
\begin{pf}
The proof is the same as the one of Lemma 6.2 in \cite{no1} replacing Lemma 10.1 by its counterpart Lemma \ref{eq:intbubbleest}.
\end{pf}
\vspace{6pt}

\noindent
Still as in \cite{no1}, the second lemma that we need for the proof of Proposition \ref{eq:expansionJ1} read as follows:
\begin{lem}\label{eq:holmos1}
Assuming the assumptions of Proposition \ref{eq:expansionJ1}, then there holds the following estimate
\begin{equation}\label{eq:esthm6}
\begin{split}
&\frac{\oint_{\partial M} K e^{3\sum_{i=1}^k\alpha_i\varphi_{a_i, \l_i}+3\sum_{r=1}^{\bar k}\beta_r v_r}wdS_g}{\oint_{\partial M} K e^{3\sum_{i=1}^k\alpha_i\varphi_{a_i, \l_i}+3\sum_{r=1}^{\bar k}\beta_r v_r}dS_g}=\\&O\left(||w||_{\mathbb{P}^{4, 3}}\left(\sum_{i=1}^k\frac{|\n_{\hat g} \mathcal{F}^{A}_i(a_i)|}{\l_i}+\sum_{i=1}^k|\alpha_i-1|\log\l_i+\sum_{r=1}^{\bar k}|\beta_r|+\sum_{i=1}^k\frac{\log \l_i}{\l_i^{2}}\right)\right).
\end{split}
\end{equation}
\end{lem}
\begin{pf}
If follows from the same arguments as in  the proof of Lemma 6.3 in \cite{no1} replacing Lemma 10.1 by  its counterpart Lemma \ref{eq:intbubbleest}.
\end{pf}
\vspace{6pt}

\noindent
Finally, as in \cite{no1}, the third and last lemma that we need for the proof of Proposition \ref{eq:expansionJ1} is the following one.
\begin{lem}\label{eq:tauiest}
Assuming the assumptions of Proposition \ref{eq:expansionJ1}, then for every \;$i=1, \cdots, k$,  there holds
\begin{equation}\label{eq:esttauieq}
\tau_i=O\left(\sum_{j=1}^k\frac{1}{\l_j}\right).
\end{equation}
\end{lem}
\begin{pf}
The proof is the same as the one Lemma 6.4 in \cite{no1} replacing Lemma 5.1 by  Lemma \ref{eq:gradientlambdaest}.
\end{pf}
\vspace{8pt}

\noindent
\begin{pfn} {of Proposition \ref{eq:expansionJ1}}\\
It follows from the same arguments as in the proof of Lemma 6.1 in \cite{no1} replacing Lemma 6.2-Lemma 6.4 in \cite{no1} by  Lemma \ref{eq:holmos}-lemma \ref{eq:tauiest} and Lemma 10.1 in \cite{no1}  by Lemma \ref{eq:intbubbleest}. Furthermore, Lemma 10.6 and Lemma 10.7 in \cite{no1} are replaced by Lemma \ref{eq:positive} and Lemma \ref{eq:positiveg}.
\end{pfn}
\vspace{6pt}

\noindent
Now, as in \cite{no1}, we have that Proposition \ref{eq:expansionJ1} implies the following direct corollaries.
\begin{cor}\label{eq:cexpansionj1}
Assuming that $\eta$ is a small positive real number with \;$0<2\eta<\varrho$\; where \;$\varrho$\; is as in \eqref{eq:cutoff}, $0<\epsilon\leq \epsilon_0$\; where \;$\epsilon_0$\; is as in \eqref{eq:mini} and \;$u:=\sum_{i=1}^k\alpha_i\varphi_{a_i, \l_i}+\sum_{r=1}^{\bar k}\beta_r(v_r-\ov{(v_r)}_{(Q, T)})$\; with the concentration  points \;$a_i$,  the masses \;$\alpha_i$, the concentrating parameters \;$\l_i$ ($i=1, \cdots, k$)  and the negativity parameters \;$\beta_r$ ($r=1, \cdots, \bar k$) satisfying \eqref{eq:afpara}, then there exists a unique \;$\bar w(\bar \alpha, A, \bar \l, \bar \beta)\in E_{A, \bar \l}$\; such that
\begin{equation}\label{eq:minj}
\mathcal{E}_g\left(u+\bar w(\bar \alpha, A, \bar \l, \bar \beta)\right)=\min_{w\in E_{A, \bar \l}, u+w\in V_R(k, \epsilon, \eta)} \mathcal{E}_g(u+w),
\end{equation}
where \;$\bar \alpha:=(\alpha_1, \cdots, \alpha_k)$, $A:=(a_1, \cdots, a_k)$, $\bar \l:=(\l_1, \cdots, \l_k)$\; and \;$\bar \beta:=(\beta_1, \cdots, \beta_k)$.\\
Furthermore, $(\bar \alpha, A, \bar \l, \bar \beta)\longrightarrow \bar w(\bar \alpha, A, \bar \l, \bar \beta)\in C^1$\; and satisfies the following estimate
\begin{equation}\label{eq:linminqua}
\frac{1}{C}||\bar w(\bar \alpha, A,  \bar \l, \bar \beta)||_{\mathbb{P}^{4, 3}}^2\leq |f\left(\bar w(\bar \alpha, A, \bar \l, \bar \beta)\right)|\leq C|| \bar w(\bar \alpha, A, \bar \l, \bar \beta)||_{\mathbb{P}^{4, 3}}^2,
\end{equation}
for some large positive constant \;$C$ independent of \;$\bar \alpha$, $A$, $\bar \l$, and \;$\bar \beta$, hence
\begin{equation}\label{eq:estbarw}
||\bar w(\bar \alpha, A, \bar \l, \bar \beta)||_{\mathbb{P}^{4, 3}}=O\left(\sum_{i=1}^k\frac{|\n_{\hat g} \mathcal{F}^{A}_i(a_i)|}{\l_i}+\sum_{i=1}^k|\alpha_i-1|\log\l_i+\sum_{r=1}^{\bar k}|\beta_r|+\sum_{i=1}^k\frac{\log \l_i}{\l_i^{2}}\right) .
\end{equation}
\end{cor}
\vspace{4pt}

\noindent
\begin{cor}\label{eq:c1expansionj1}
Assuming that \;$\eta$\; is a small positive real number with \;$0<2\eta<\varrho$\; where \;$\varrho$\; is as in \eqref{eq:cutoff}, $0<\epsilon\leq \epsilon_0$\; where \;$\epsilon_0$\; is as in \eqref{eq:mini}, and \;$u_0:=\sum_{i=1}^k\alpha_i^0\varphi_{a_i^0, \l_i^0}+\sum_{r=1}^{\bar k}\beta_r^0(v_r-\ov{(v_r)}_{(Q, T)})$\; with the concentration  points \;$a_i^0$,  the masses \;$\alpha_i^0$, the concentrating parameters \;$\l_i^0$ ($i=1, \cdots, k$) and the negativity parameters \;$\beta_r^0$ ($r=1, \cdots, \bar k$) satisfying \eqref{eq:afpara}, then there exists an open neighborhood  \;$U$\; of $(\bar\alpha^0, A^0, \bar \l^0, \bar \beta^0)$  (with \;$\bar \alpha^0:=(\alpha^0_1, \cdots, \alpha^0_k)$, $A^0:=(a_1^0, \cdots, a^0_k)$,  $\bar \l:=(\l_1^0, \cdots, \l_k^0)$\; and  \;$\bar \beta^0:=(\beta_1^0, \cdots, \beta_{\bar k}^0)$) such that for every \;$(\bar \alpha, A, \bar \l, \bar \beta)\in U$\; with \;$\bar \alpha:=(\alpha_1, \cdots, \alpha_k)$, $A:=(a_1, \cdots, a_k)$,  $\bar \l:=(\l_1, \cdots, \l_k)$, $\bar \beta:=(\beta_1, \cdots, \beta_{\bar k})$, and the \;$a_i$,  the \;$\alpha_i$, the \;$\l_i$ ($i=1, \cdots, k$)  and the \;$\beta_r$ ($r=1, \cdots, \bar k$) satisfying \eqref{eq:afpara}, and \;$w$\; satisfying \eqref{eq:afpara}  with \;$\sum_{i=1}^k\alpha_i\varphi_{a_i, \l_i}+\sum_{r=1}^{\bar k}\beta_r(v_r-\ov{(v_r)}_{(Q, T)})+w\in V_R(k, \epsilon, \eta)$, we have the existence of a change of variable
\begin{equation}\label{eq:changev}
w\longrightarrow V
\end{equation}
from a neighborhood of \;$\bar w(\bar \alpha, A, \bar \l, \bar \beta)$\; to a neighborhood of \;$0$\; such that
\begin{equation}\label{eq:expjv}
\begin{split}
&\mathcal{E}_g\left(\sum_{i=1}^k\alpha_i\varphi_{a_i, \l_i}+\sum_{r=1}^{\bar k}\beta_r(v_r-\ov{(v_r)}_{(Q, T)})+w\right)=\\&\mathcal{E}_g\left(\sum_{i=1}^k\alpha_i\varphi_{a_i, \l_i}+\sum_{r=1}^{\bar k}\beta_r(v_r-\ov{(v_r)}_{(Q, T)})+\bar w (\bar \alpha, A, \bar \l, \bar \beta)\right)\\&+\frac{1}{2}\partial^2 \mathcal{E}_g\left(\sum_{i=1}^k\alpha_i^0\varphi_{a_i^0, \l_i^0}+\sum_{r=1}^{\bar k}\beta_r^0(v_r-\ov{(v_r)}_{(Q, T)})+\bar w(\bar \alpha^0, A^0, \bar \l^0, \bar \beta^0)\right)(V, V),
\end{split}
\end{equation}
\end{cor}
\vspace{6pt}

\noindent
Thus, as in \cite{no1}, with this new variable, it is easy to see that in \;$V_R(k, \epsilon, \eta)$\; we have a splitting of the variables \;$(\bar \alpha, A, \bar \l, \bar \beta)$\; and \;$V$, namely that one can decrease the Euler-Lagrange functional \;$\mathcal{E}_g$\; in the variable \;$V$\; without touching the variable \;$(\bar \alpha, A, \bar \l, \bar \beta)$\; by considering just the flow
\begin{equation}\label{eq:vflow}
\frac{dV}{dt}=-V.
\end{equation}
So, as in \cite{no1}, and for the same reasons, to develop a Morse theory for \;$\mathcal{E}_g$\;  is equivalent to do one for the functional 
\begin{equation}\label{eq:finitedj}
\bar{\mathcal{E}_g}(\bar \alpha, A, \bar \l, \bar \beta):=\mathcal{E}_g\left(\sum_{i=1}^k\alpha_i\varphi_{a_i, \l_i}+\sum_{r=1}^{\bar k}\beta_r(v_r-\ov{(v_r)}_{(Q, T)})+\bar w(\bar \alpha, A, \bar \l, \bar \beta)\right),
\end{equation}
where \;$\bar \alpha=(\alpha_1, \cdots, \alpha_k)$, $A=(a_1, \cdots, a_k)$, $\bar \l=(\l_1, \cdots, \l_k)$\; and \;$\bar \beta=\beta_1, \cdots, \beta_{\bar k}$\; with the concentration  points \;$a_i$,  the masses \;$\alpha_i$, the concentrating parameters \;$\l_i$ ($i=1, \cdots, k$)  and the negativity parameters \;$\beta_r$ ($r=1, \cdots, \bar k$) satisfying \eqref{eq:afpara}, and \;$\bar w(\bar \alpha, A, \bar \l, \bar \beta)$\; is as in Corollary \ref{eq:cexpansionj1}. 
\vspace{6pt}

\noindent
Finally, we have the following energy estimate of \;$\mathcal{E}_g$\; on \;$V_R(k, \epsilon, \eta)$ .
\begin{lem}\label{eq:expansionj}
Under the assumptions of Proposition \ref{eq:expansionJ1}, $\forall u=\ov{u}_{(Q, T)}+\sum_{i=1}^k\alpha_i\varphi_{a_i, \l_i}+\sum_{r=1}^{\bar k}\beta_r (v_r-\ov{(v_r)}_{(Q, T)})+w\in V_R(k, \epsilon, \eta)$, we have 
\begin{equation*}
\begin{split}
&\mathcal{E}_g(u)=\mathcal{E}_g\left(\sum_{i=1}^k\alpha_i\varphi_{a_i, \l_i}+\sum_{r=^1}^{\bar k}\beta_r(v_r-\ov{(v_r)}_{(Q, T)})+w\right)=C^{k}_0-8\pi^2\mathcal{F}_K(a_1, \dots, a_k)+2\sum_{r=1}^{\bar k}\mu_r\beta_r^2\\&+\sum_{i=1}^k(\alpha_i-1)^2\left[16\pi^2\log\l_i+8\pi^2H(a_i, a_i)+C^{k}_1\right]-\frac{c^18\pi^2}{9}\sum_{i=1}^k\frac{1}{\l_i^2}\left(\frac{\D_{\hat g_{a_i}}\mathcal{F}^{A}_i(a_i)}{\mathcal{F}^{A}_i(a_i)}-\frac{3}{4}R_{\hat g}(a_i)\right)\\&+\frac{1}{2}\partial^2 \mathcal{E}_g\left(\sum_{i=1}^k\alpha_i^0\varphi_{a_i^0, \l_i^0}+\sum_{r=1}^{\bar k}\beta_r^0(v_r-\ov{(v_r)}_{(Q, T)})+\bar w(\bar \alpha^0, A^0, \bar \l^0, \bar \beta^0)\right)(V, V)\\&+8\pi^2\sum_{i=1}^k(\alpha_i-1)\left[\sum_{r=1}^{\bar k}2\beta_r\left(v_r-\ov{(v_r)}_{(Q, T)}\right)(a_i)-\sum_{ j=1, j\neq i}^k (\alpha_j-1)G(a_i, a_j)\right]\\&+\frac{c^18\pi^2}{9}\sum_{i=1}^k\frac{\tilde \tau_i}{\l_i^2}\left(\frac{\D_{\hat g_{a_i}}\mathcal{F}^{A}_i(a_i)}{\mathcal{F}^{A}_i(a_i)}-\frac{3}{4}R_{\hat g}(a_i)\right)\\&+\frac{16\pi^2}{3}\sum_{i=1}^k\log(1-\tilde\tau_i)+O\left(\sum_{i=1}^k|\alpha_i-1|^3+\sum_{r=1}^{\bar k}|\beta_r|^3+\sum_{i=1}^k\frac{1}{\l^3_i}
+||\bar w(\bar \alpha, A, \bar \l, \bar \beta)||^2_{\mathbb{P}^{4, 3}}\right),
\end{split}
\end{equation*}
where \;$O\left(1\right)$ means here \;$O_{\bar\alpha, A, \bar \l, \bar \beta, \epsilon}\left(1\right)$ \;with \;$\bar\alpha=(\alpha_1, \cdots, \alpha_k)$, \;$A:=(a_1, \cdots, a_k)$, \;$\bar \l:=(\l_1, \cdots, \l_k)$, \;$\bar \beta:=(\beta_1, \cdots, \beta_{\bar k})$\; and for \;$i=1, \cdots, k$, \;$\tilde \tau_i$\; is as in Lemma \ref{eq:energyest}.
where \;$\bar w(\bar \alpha, A, \bar \l, \bar \beta)$\; is as in Corollary \ref{eq:cexpansionj1}.
\end{lem}
\begin{pf}
It follows directly from Lemma \ref{eq:energyest}, formula \eqref{eq:expjv} and Proposition \ref{eq:expansionJ1}.
\end{pf}
\subsection{Morse lemma at infinity}
In this subsection, we derive a Morse Lemma at infinity for \;$\mathcal{E}_g$. As in \cite{no1}, in order to do that, we first construct a pseudo-gradient for \;$\bar{\mathcal{E}_g}(\bar \alpha, A, \bar \l, \bar \beta)$, where \;$\bar{\mathcal{E}_g}(\bar \alpha, A, \bar \l, \bar \beta)$\; is defined as in \eqref{eq:finitedj} exploiting the gradient estimates derived previously. Indeed, we have:
\begin{pro}\label{eq:conspseudograd}
Assuming that $\eta$ is a small positive real number with \;$0<2\eta<\varrho$\; where \;$\varrho$\; is as in \eqref{eq:cutoff}, and \;$0<\epsilon\leq \epsilon_0$\; where \;$\epsilon_0$\; is as in \eqref{eq:mini}, then there exists a pseudogradient \;$\mathcal{W}_g$\; of \;$\bar{\mathcal{E}_g}(\bar \alpha, A, \bar \l, \bar \beta)$\; such that \\
1) For every \;$u:=\sum_{i=1}^k\alpha_i\varphi_{a_i, \l_i}+\sum_{r=1}^{\bar k}\beta_r (v_r-\ov{(v_r)}_{(Q, T)})\in V_R(k, \epsilon, \eta)$\; with the concentration  points \;$a_i$, the masses \;$\alpha_i$, the concentrating parameters \;$\l_i$ ($i=1, \cdots, k$)  and the negativity parameters \;$\beta_r$ ($r=1, \cdots, \bar k$)\; satisfying \eqref{eq:afpara}, there holds
\begin{equation}\label{eq:pseudoexact}
\large\left\langle-\n^{\mathbb{P}^{4, 3}} \mathcal{E}_g(u), \;\mathcal{W}_g\large\right\rangle_{\mathbb{P}^{4, 3}}\geq c\left(\sum_{i=1}^k\frac{1}{\l_i^2}+\sum_{i=1}^k\frac{|\n_{\hat g}\mathcal{F}^{A}_i(a_i)|}{\l_i}+\sum_{i=1}^k|\alpha_i-1|+\sum_{i=1}^k|\tau_i|+\sum_{r=1}^{\bar k}|\beta_r|)\right),
\end{equation}
and  for every \;$u:=\sum_{i=1}^k\alpha_i\varphi_{a_i, \l_i}+\sum_{r=1}^{\bar k}\beta_r \l(v_r-\ov{(v_r)}_{(Q, T)})+\bar w(\bar \alpha, A, \bar \l, \bar \beta)\in V_R(k, \epsilon, \eta)$ with the concentration  points \;$a_i$, the masses \;$\alpha_i$, the concentrating parameters \;$\l_i$ ($i=1, \cdots, k$)  and the negativity parameters \;$\beta_r$ ($r=1, \cdots, \bar k$) satisfying \eqref{eq:afpara}, and \;$\bar w(\bar \alpha, A, \bar \l, \bar \beta)$\; is as in \eqref{eq:minj},  there holds
\begin{equation}\label{eq:pseudoperturb}
\left\langle-\n^{\mathbb{P}^{4, 3}} \mathcal{E}_g(u+\bar w), \;\mathcal{W}_g+\frac{\partial \bar w}{\partial (\bar\alpha, A, \bar \l, \bar \beta)}\right\rangle_{_{\mathbb{P}^{4, 3}}}\geq c\left(\sum_{i=1}^k\frac{1}{\l_i^2}+\sum_{i=1}^k\frac{|\n_{\hat g}\mathcal{F}^{A}_i(a_i)|}{\l_i}+\sum_{i=1}^k|\alpha_i-1|+\sum_{i=1}^k|\tau_i|+\sum_{r=1}^{\bar k}|\beta_r|\right),
\end{equation}
where \;$c$\; is a small positive constant independent of \;$A:=(a_1, \cdots, a_k)$, $\bar\alpha=(\alpha_1, \cdots, \alpha_k)$, $\bar\l=(\l_1, \cdots, \l_k)$, $\bar \beta=(\beta_1, \cdots, \beta_{\bar k})$\; and \;$\epsilon$.\\
2) $\mathcal{W}_g$\; is a \;$||\cdot||_{\mathbb{P}^{4, 3}}$-bounded vector field and is compactifying outside the region where \;$A$\; is very close to a critical point \;$B$\; of \;$\mathcal{F}_K$\; satisfying \;$\mathcal{L}_K(B)<0$.
\end{pro}
\begin{pf}
It follows from the same arguments as in the proof of Proposition 8.1 in \cite{no1} replacing formulas (52)-(54), Lemma 5.1, Corollary 5.2 and Lemma 5.3-Lemma 5.5 in \cite{no1}  with \eqref{eq:relationderivative}-\eqref{eq:auxiindexa1}, Lemma \ref{eq:gradientlambdaest},  Corollary \ref{eq:cgradientlambdaest} and Lemma \ref{eq:gradientalpha}-Lemma \ref{eq:gradientbeta}. Furthermore,  Lemma 4.1, Lemma 7.1, and Lemma 10.5 in \cite{no1} are replaced by Lemma \ref{eq:energyest}, Lemma \ref{eq:expansionJ1} and Lemma \ref{eq:auxibubbleest1}
\end{pf}

\vspace{6pt}

 \noindent
Now, as in \cite{no1}, we have that Proposition \ref{eq:conspseudograd} implies the following characterization of the  critical points at infinity of \;$\mathcal{E}_g$.
 \begin{cor}\label{eq:loccritinf}
 1) The  critical points at infinity of \;$\mathcal{E}_g$\; correspond to the "configurations" \;$\alpha_i=1$, $\l_i=+\infty$, \;$\tau_i=0$ \;$i=1, \cdots, k$, \;$\beta_r=0$, \;$r=1, \cdots, \bar k$, \;$A$\; is a critical point of \;$\mathcal{F}_K$\; and  \;$V=0$, and we denote them by \;$z^{\infty}$\; with \;$z$\; being the corresponding critical point of $\mathcal{F}_K$ .\\\\
 2) The ``true'' critical points at infinity of \;$\mathcal{E}_g$\; are the \;$z^{\infty}$\; satisfying \;$\mathcal{L}_K(z)<0$\; and we denote them by \;$x^{\infty}$\; with \;$x$\; being the corresponding critical point of \;$\mathcal{F}_K$.\\\\
 3) The ``false'' critical points at infinity of \;$\mathcal{E}_g$\; are the \;$z^{\infty}$\; satisfying \;$\mathcal{L}_K(z)>0$\; and we denote them by \;$y^{\infty}$\; with \;$y$\; being the corresponding critical point of \;$\mathcal{F}_K$.\\ 
 4) The \;$\mathcal{E}_g$-energy of a critical point at infinity \;$z^{\infty}$\; denoted by \;$\mathcal{J}_{g}(z^{\infty})$\; is given by
 \begin{equation}\label{eq:infinitycriticallevel}
 \mathcal{J}_{g}(z^{\infty})=C_0^{k}-8\pi^2\mathcal{F}_K(z_1, \dots, z_k)
 \end{equation}
 where \;$z=(z_1, \cdots, z_k)$\; and \;$C_0^{k}$\; is as in Lemma \ref{eq:energyest}.
 \end{cor}
 \begin{pf}
 Point 1)- Point 3) follow from \eqref{eq:auxiindexa1}, Lemma \ref{eq:deformlemr}, the discussions right after \eqref{eq:expjv}, and Proposition \ref{eq:conspseudograd}, while Point 4) follows from Point 1) combined with \eqref{eq:estbarw} and Lemma \ref{eq:expansionj}. 
\end{pf}
 \vspace{6pt}

 \noindent
Finally, we are going to conclude this subsection by establishing an analogue of the classical Morse lemma for both ``true'' and ``false'' critical points at infinity. In order to do that, we first remark that, as in \cite{no1}, the arguments of Proposition \ref{eq:conspseudograd} implies that \;$V_{-}:=\{u\in V_R(k, \epsilon, \eta): \; l_K(A)<0, \;\forall  r\in \{1, \cdots, \bar k\}\,\,\;\,|\beta_r|\leq 2\tilde C_0\left(\sum_{i=1}^k\frac{|\n_{\hat g} \mathcal{F}^{A}(a_i)|}{\l_i} +\sum_{i=1}^k|\alpha_i-1|+\sum_{i=1}^k|\tau_i|+\sum_{i=1}^k\frac{1}{\l_i^2}\right),\;\;\forall \;i\in \{1, \cdots, k\}\;\,\;|\tau_i|\leq 2\frac{\hat C_0}{\l_i^2},  \;\; \text{and}\;\;\;\forall i\in \{1, \cdots, k\}\;\;\frac{|\n_{\hat g} \mathcal{F}^{A}_i(a_i)|}{\l_i}\leq 4\frac{C_0}{\l_i^2}\;\} $ and $V_{+}:=\{u\in V_R(k, \epsilon, \eta): l_K(A)>0, \forall  r\in \{1, \cdots, \bar k\}\,\,\;\,|\beta_r|\leq 2\tilde C_0\left(\sum_{i=1}^k\frac{|\n_{\hat g} \mathcal{F}^{A}(a_i)|}{\l_i} +\sum_{i=1}^k|\alpha_i-1|+\sum_{i=1}^k|\tau_i|+\sum_{i=1}^k\frac{1}{\l_i^2}\right),\\\;\;\forall \;i\in \{1, \cdots, k\}\;\,\;|\tau_i|\leq 2\frac{\hat C_0}{\l_i^2},  \;\; \text{and}\;\;\;\forall i\in \{1, \cdots, k\}\;\;\frac{|\n_{\hat g} \mathcal{F}^{A}_i(a_i)|}{\l_i}\leq 4\frac{C_0}{\l_i^2}\;\} $ (where $\tilde C_0$, $\hat C_0$ and $C_0$ are large positive constants) are respectively a neighborhood of the ``true''  and ``false'' critical points at infinity of the variational problem. 
Hence, as in \cite {no1}, \eqref{eq:estbarw}, Corollary\ref{eq:c1expansionj1}, Lemma \ref{eq:expansionj} and classical Morse lemma imply the following Morse type lemma for a ``true'' critical point at infinity.
 \begin{lem}\label{eq:morleminf}(Morse lemma at infinity near a ``true'' one)\\
Assuming that $\eta$ is a small positive real number with \;$0<2\eta<\varrho$\; where \;$\varrho$\; is as in \eqref{eq:cutoff}, \;$0<\epsilon\leq \epsilon_0$\; where \;$\epsilon_0$\; is as in \eqref{eq:mini} and \;$u_0:=\sum_{i=1}^k\alpha_i^0\varphi_{a_i^0, \l_i^0}+\sum_{r=1}^{\bar k}\beta_r^0(v_r-\ov{(v_r)}_{(Q, T)})+\bar w((\bar\alpha^0, A^0, \bar \l^0, \bar \beta^0))\in V_{-}(k, \epsilon, \eta)$ (where \;$\bar \alpha^0:=(\alpha^0_1, \cdots, \alpha^0_k)$, $A^0:=(a_1^0, \cdots, a^0_k)$,  $\bar \l:=(\l_1^0, \cdots, \l_k^0)$ and  $\bar \beta^0:=(a_1^0, \cdots, \beta_{\bar k}^0)$) with the concentration  points \;$a_i^0$,  the masses\; $\alpha_i^0$, the concentrating parameters \;$\l_i^0$ ($i=1, \cdots, k$) and the negativity parameters \;$\beta_r^0$ ($r=1, \cdots, \bar k$) satisfying \eqref{eq:afpara} and furthermore \;$A^0\in Crit(\mathcal{F}_K)$, then there exists an open neighborhood  \;$U$\; of \;$(\bar\alpha^0, A^0, \bar \l^0, \bar \beta^0)$\; such that for every \;$(\bar \alpha, A, \bar \l, \bar \beta)\in U$\; with \;$\bar \alpha:=(\alpha_1, \cdots, \alpha_k)$, $A:=(a_1, \cdots, a_k)$,  $\bar \l:=(\l_1, \cdots, \l_k)$, $\bar \beta:=(\beta_1, \cdots, \beta_{\bar k})$, and the \;$a_i$,  the \;$\alpha_i$, the \;$\l_i$ ($i=1, \cdots, k$)  and the \;$\beta_r$ ($r=1, \cdots, \bar k$) satisfying \eqref{eq:afpara}, and \;$w$\; satisfying \eqref{eq:afpara}  with \;$u=\ov{u}_{(Q, T)}\sum_{i=1}^k\alpha_i\varphi_{a_i, \l_i}+\sum_{r=1}^{\bar k}\beta_r(v_r-\ov{(v_r)}_{(Q, T)})+w\in V_{-}(k, \epsilon, \eta)$, we have the existence of a change of variable
\begin{equation}\label{eq:morsevinf}
\begin{split}
&\alpha_i\longrightarrow s_i ,i=1, \cdots, k,\\& A\longrightarrow \tilde A=(\tilde A_{-}, \tilde A_{+})\\&\l_1\longrightarrow \theta_1,\\&\tau_i\longrightarrow \theta_i, i=2, \cdots, k,\\&\beta_r\longrightarrow \tilde \beta_r\\ &V\longrightarrow \tilde V,
\end{split}
\end{equation}
such that
\begin{equation}
\mathcal{E}_g(u)=\mathcal{E}_g\left(\sum_{i=1}^k\alpha_i\varphi_{a_i, \l_i}+\sum_{r=1}^{\bar k}\beta_r\left(v_r-\ov{(v_r)}_{(Q, T)}\right)+w\right)=-|\tilde A_{-}|^2+|\tilde A_{+}|^2+\sum_{i=1}^{k}s_i^2-\sum_{r=1}^{\bar k}\tilde \beta_r^2+\theta_1^2-\sum_{i=2}^k\theta_i^2+||\tilde V||^2
\end{equation}
where \;$\tilde A=(\tilde A_{-}, \;\tilde A_{+})$\; is the Morse variable of the map \;$\mathcal{J}_{g}: (\partial M)^k\setminus F((\partial M)^k)\longrightarrow \R$\; which is defined by the right hand side of \eqref{eq:infinitycriticallevel}. Hence a ``true``critical point at infinity \;$x^{\infty}$\; of \;$\mathcal{E}_g$\; has Morse index at infinity  \;$\mathcal{M}_{\infty}(x^{\infty})=i_{\infty}(x)+\bar k$.
\end{lem}
\vspace{6pt}

\noindent
Similarly, and for the same reasons as above, we have the following analogue of the classical Morse lemma for a ''false`` critical point at infinity.
\begin{lem}\label{eq:morleminf1}(Morse lemma at infinity near a ``false'' one)\\
Assuming that \;$\eta$\; is a small positive real number with \;$0<2\eta<\varrho$\; where \;$\varrho$\; is as in \eqref{eq:cutoff}, $0<\epsilon\leq \epsilon_0$\; where \;$\epsilon_0$\; is as in \eqref{eq:mini} and \;$u_0:=\sum_{i=1}^k\alpha_i^0\varphi_{a_i^0, \l_i^0}+\sum_{r=1}^{\bar k}\beta_r^0(v_r-\ov{(v_r)}_{(Q, T)})+\bar w((\bar\alpha^0, A^0, \bar \l^0, \bar \beta^0))\in  V_{+}(k, \epsilon, \eta) $ (where \;$\bar \alpha^0:=(\alpha^0_1, \cdots, \alpha^0_k)$, \;$A^0:=(a_1^0, \cdots, a^0_k)$,  \;$\bar \l:=(\l_1^0, \cdots, \l_k^0)$\; and  \;$\bar \beta^0:=(a_1^0, \cdots, \beta_{\bar k}^0)$) with the concentration  points \;$a_i^0$,  the masses \;$\alpha_i^0$, the concentrating parameters \;$\l_i^0$ ($i=1, \cdots, k$) and the negativity parameters \;$\beta_r^0$ ($r=1, \cdots, \bar k$) satisfying \eqref{eq:afpara} and furthermore \;$A^0\in Crit(\mathcal{F}_K)$, then there exists an open neighborhood  \;$U$\; of $(\bar\alpha^0, A^0, \bar \l^0, \bar \beta^0)$ such that for every \;$(\bar \alpha, A, \bar \l, \bar \beta)\in U$\; with \;$\bar \alpha:=(\alpha_1, \cdots, \alpha_k)$, $A:=(a_1, \cdots, a_k)$,  $\bar \l:=(\l_1, \cdots, \l_k)$, $\bar \beta:=(\beta_1, \cdots, \beta_{\bar k})$, and the \;$a_i$,  the \;$\alpha_i$, the $\l_i$ ($i=1, \cdots, k$)  and the \;$\beta_r$ ($r=1, \cdots, \bar k$) satisfying \eqref{eq:afpara}, and \;$w$\; satisfying \eqref{eq:afpara}  with\; $u=\ov{u}_{(Q, T)}+\sum_{i=1}^k\alpha_i\varphi_{a_i, \l_i}+\sum_{r=1}^{\bar k}\beta_r(v_r-\ov{(v_r)}_{Q^n})+w\in V_{+}(k, \epsilon, \eta)$, we have the existence of a change of variable
\begin{equation}\label{eq:morsevinf1}
\begin{split}
&\alpha_i\longrightarrow s_i ,i=1, \cdots, k,\\& A\longrightarrow \tilde A=(\tilde A_{-}, \tilde A_{+})\\&\l_1\longrightarrow \theta_1,\\&\tau_i\longrightarrow \theta_i, i=2, \cdots, k,\\&\beta_r\longrightarrow \tilde \beta_r\\ &V\longrightarrow \tilde V,
\end{split}
\end{equation}
such that
\begin{equation}
\mathcal{E}_g(u)=\mathcal{E}_g\left(\sum_{i=1}^k\alpha_i\varphi_{a_i, \l_i}+\sum_{r=1}^{\bar m}\beta_r\left(v_r-\ov{(v_r)}_{(Q, T)}\right)+w\right)=-|\tilde A_{-}|^2+|\tilde A_{+}|^2+\sum_{i=1}^{k}s_i^2-\sum_{r=1}^{\bar k}\tilde \beta_r^2-\sum_{i=1}^k\theta_i^2+||\tilde V||^2,
\end{equation}
where \;$\tilde A=(\tilde A_{-}, \tilde A_{+})$\; is the Morse variable of the map \;$\mathcal{J}_{g}: (\partial M)^k\setminus F((\partial M)^k)\longrightarrow \R$\; which is defined by the right hand side of \eqref{eq:infinitycriticallevel}. Hence a ``false``critical point at infinity \;$y^{\infty}$\; of \;$\mathcal{E}_g$\; has Morse index at infinity  \;$M_{\infty}(y^{\infty})=i_{\infty}(y)+1+\bar k$.
\end{lem}

\section{Proof of existence theorems}
In this section, we show how the Morse lemma at infinity implies the main existence results via strong Morse type inequalities or Barycenter technique of Bahri-Coron.
\subsection{Topology of vey high and negative sublevels of \;$\mathcal{E}_g$} 
We study the topology of very high sublevels of \;$\mathcal{E}_g$\; and its every negative ones.  We start with the very high sublevels of \;$\mathcal{E}_g$\; and first derive the following Lemma. 
\begin{lem}\label{eq:energybddinf}
Assuming that \;$\eta$\; is a small positive real number with \;$0<2\eta<\varrho$\; where \;$\varrho$\; is as in \eqref{eq:cutoff},  then there exists \;$\hat C_0^k:=\hat C^k_0(\eta)$ such that  for every  \,$0<\epsilon\leq \epsilon_0$\; where \;$\epsilon_0$\; is as in \eqref{eq:mini}, there holds
\begin{equation*}
V(k, \epsilon, \eta)\subset (\mathcal{E}_g)^{\hat C_0^k}-(\mathcal{E}_g)^{-\hat C_0^k}.
\end{equation*}
\end{lem}
\begin{pf}
It follows directly from \eqref{eq:para}-\eqref{eq:afpara}, Proposition \ref{eq:expansionJ1}, Lemma \ref{eq:tauiest}, and Lemma \ref{eq:expansionj}.
\end{pf}
\vspace{4pt}

\noindent
Next, combining Proposition \ref{eq:escape} and the latter lemma, we have the following corollary.
\begin{cor}\label{eq:energybddcrit}
There exists a large positive constant \;$\hat C_1^k$\; such that
\begin{equation*}
Crit(\mathcal{E}_g)\subset (\mathcal{E}_g)^{\hat C_1^k}-(\mathcal{E}_g)^{-\hat C_1^k}.
\end{equation*}
\end{cor}
\begin{pf}
It follows, via a contradiction argument, from the the fact that \;$\mathcal{E}_g$\; is invariant by translation by constants, Proposition \ref{eq:escape}, and  Lemma \ref{eq:energybddinf}.
\end{pf}
\vspace{6pt}

\noindent
Now, we are ready to characterize the topology of very high sublevels of \;$\mathcal{E}_g$. Indeed, as in \cite{no1} and for the same reasons, we have that Lemma \ref{eq:deformlemr}, Lemma \ref{eq:energybddinf} and Corollary \ref{eq:energybddcrit} imply the following one which describes the topology of very high sublevels of the Euler-Lagrange functional \;$\mathcal{E}_g$. 
\begin{lem}\label{eq:tophigh}
Assuming that \;$\eta$\; is a small positive real number with \;$0<2\eta<\varrho$\; where \;$\varrho$\; is as in \eqref{eq:cutoff}, then there exists a large positive constant  \;$L^k:=L^k(\eta)$ with \;$L^k>2\max\{\hat C_0^k, \hat C_1^k\}$\;  such that for every \;$L\geq L^k$, we have  that \;$(\mathcal{E}_g)^L$\; is a deformation retract of \;$\mathcal{H}_{\frac{\partial }{\partial n}}$, and hence it has the homology of a point, where \;$\hat C^k_0$\; is as in Lemma \ref{eq:energybddinf} and \;$\hat C^k_1$\; is as in Lemma \ref{eq:energybddcrit}.
\end{lem}
\vspace{6pt}

\noindent
Next, we turn to the study of the topology of very negative sublevels of \;$\mathcal{E}_g$\; when \;$k\geq 2$\; or \;$\bar k\geq 1$. Indeed, as in \cite{no1} and for the same reasons, we have that the well-know topology of very negative sublevels in the {\em nonresont} case (see \cite{nd6}), Proposition \ref{eq:escape}, Lemma \ref{eq:energybddinf} and Corollary \ref{eq:energybddcrit} imply the following lemma which gives the homotopy type of the very negative sublevels of the Euler-Lagrange functional \;$\mathcal{E}_g$.
\begin{lem}\label{eq:topnegative}
Assuming that \;$k\geq 2$\; or \;$\bar k\geq 1$, and \;$\eta$\; is a small positive real number with \;$0<2\eta<\varrho$\; where \;$\varrho$\; is as in \eqref{eq:cutoff}, then there exists a large positive constant \;$L_{k, \bar k}:=L_{k, \bar k}(\eta)$\; with  \;$L_{k, \bar k}>2\max\{\hat C_0^k, \hat C_1^k\}$\; such that for every \;$L\geq L_{k, \bar k}$, we have  that \;$(\mathcal{E}_g)^{-L}$\;  has the same homotopy type as  \;$B_{k-1}(\partial M)$\; if \;$k\geq 2$\; and \;$\bar k=0$,  as\;$A_{k-1, \bar k}$\; if \;$k\geq 2$\; and \;$\bar k\geq 1$\; and  as \;$S^{\bar k-1}$\; if \;$k=1$\; and \;$\bar k\geq 1$, where \;$\hat C^k_0$\; is as in Lemma \ref{eq:energybddinf} and \;$\hat C^k_1$\; as in Lemma \ref{eq:energybddcrit} .
\end{lem}
\vspace{6pt}

\noindent
However, as in \cite{nd4}, to prove Theorem \ref{eq:existence}, we need a further information about the topology of very negative sublevels of \;$\mathcal{E}_g$. In order to derive that, we first make some definitions. For \;$p\in \N^*$\; and \;$\lambda>0$, we define $$f_p(\lambda) : B_p(\partial M)\longrightarrow \mathcal{H}_{\frac{\partial }{\partial n}}$$ as follows
\begin{equation}\label{eq:fplambda}
f_p(\lambda)(\sum_{i=1}^p \alpha_i\delta_{a_i}):=\sum_{i=1}^p \alpha_i \varphi_{a_i, \lambda},\;\;\;\;\sigma=\sum_{i=1}^p \alpha_i\delta_{a_i}\in B_p(\partial M),
\end{equation}
with the \;$\varphi_{a_i, \l}$'s defined by \eqref{eq:projbubble}. Furthermore, when  \;$\bar k\geq 1$, for  \;$\Theta>0$, we define
\begin{equation}
\Psi_{p, \bar k}(\lambda, \Theta): A_{p, \bar k}\longrightarrow \mathcal{H}_{\frac{\partial }{\partial n}}
\end{equation}
as follows
\begin{equation}
\Psi_{p, \bar k}(\lambda, \theta)(\sigma, s):=
\begin{cases}
\varphi_s+f_{p}(\lambda)(\sigma)\;\;\;&\text{for}\;\;|s|\leq \frac{1}{4}, \;\;\;\sigma\in B_p(\partial M),\\
\varphi_s+f_p(2\lambda-1+4(1-\lambda)|s|)(\sigma)\;\;\;&\text{for} \;\;\frac{1}{4}\leq |s|\leq\frac{1}{2}, \;\;\;\sigma\in B_p(\partial M),\\
\varphi_s+2(1-f_p(1)(\sigma))|s|+2 f_p(1)-1\;\;\;&\text{for}\;\;\;|s|\geq \frac{1}{2},\;\;\;\sigma\in B_p(\partial M),
\end{cases}
\end{equation}
where \;$\varphi_s$\; is defined by the following formula
\begin{equation}\label{eq:varphis}
\varphi_s=\Theta\sum_{r=1}^{\bar k}s_r (v_r-\ov{(v_r)}_{(Q, T)}),
\end{equation}
with \;$s=(s_1,\cdots,s_{\bar k})$. As in \cite{nd4}, concerning the \;$f_p(\l)$'s, we have the following estimates.
\begin{lem}\label{eq:mappingbary}
Assuming that \;$p\in \N^*$, then we have\\
1) If \;$p<k$, then for every \;$L>0$, there exists \;$\lambda^L_p>0$\; such that for all \;$\lambda\ge\lambda^L_p$, we have
$$
f_p(\lambda)(B_p(\partial M))\subset (\mathcal{E}_g)^{-L}.
$$
2) If \;$p=k$, then there exist \;$\hat C_k>0$\; and \;$\lambda_k>0$\; such that for all \;$\lambda\ge\lambda_k$, we have
$$
f_k(\lambda)(B_k(\partial M))\subset (\mathcal{E}_g)^{\hat C_k}.
$$
3) There exists \;$\hat C^k>0$\; such that up to taking $\;\epsilon_0$\; smaller, where \;$\epsilon_0$\; is given by \eqref{eq:mini}, we have that for every \;$0<\epsilon\leq \epsilon_0$, there holds
$$
V(k, \epsilon)\subset (\mathcal{E}_g)^{\hat C^k}.
$$
\end{lem}
\begin{pf}
It follows from the same arguments as in the proof of Lemma 3.1 in \cite{nd4} by using Lemma \ref{eq:intbubbleest}, Lemma \ref{eq:selfintest} and Lemma \ref{l:A1}-Lemma \ref{eq:lemma3}.
\end{pf}
\vspace{6pt}

\noindent
Still, as in \cite{nd4}, we have the following estimates for the \;$\Psi_p(\l, \Theta)$'s when \;$\bar k\geq 1$.
\begin{lem}\label{eq:mappingbarysphere}
Assuming that $p\in \N^*$, then we have\\
1) If\; $1\leq p<k$, then for every \;$L>0$, there exists \;$\lambda^L_{p, \bar k}>0$\;  and \;$\Theta_{p, \bar k}^L>0$\; such that for all \;$\lambda\geq\lambda^L_{p, \bar k}$, we have
$$
\Psi_{p, \bar k}(\lambda, \Theta_{p, \bar k}^L)(A_{p, \bar k})\subset (\mathcal{E}_g)^{-L}.
$$
2) If \;$p=k$\; and \;$\Theta>0$, then there exists \;$C_{k,\bar k}^{\Theta}>0$, \;$\l_{k, \bar k}^{\Theta}>0$, such that for every \;$\l\ge\l^{\Theta}_{k, \bar k}$, we have
$$
\Psi_{k, \bar k}(\lambda, \Theta)(A_{k, \bar k})\subset (\mathcal{E}_g)^{C_{k, \bar k}^{\Theta}}.
$$
3) If \;$\Theta>0$ , then there exists \;$C^{k, \bar k}_{\Theta}>0$\; such that up to taking \;$\epsilon_0$\; smaller, where \;$\epsilon_0$\; is given by \eqref{eq:mini}, we have that for every \;$0<\epsilon\leq \epsilon_0$, there holds
$$
V(k, \epsilon, \Theta)\subset (\mathcal{E}_g)^{C^{k, \bar k}_{\Theta}}.
$$
\end{lem}
\begin{pf}
It follows from the same arguments as in the proof of Lemma 4.1 in \cite{nd4} by using Lemma \ref{eq:mappingbary}.
\end{pf}
\vspace{6pt}

\noindent
On the other hand, as in \cite{nd4}, Lemma \ref{eq:topnegative} and Lemma \ref{eq:mappingbary} imply the following one:
\begin{lem}\label{eq:verynegative}
Assuming that \;$k\geq 2$, $\bar k=0$, and  \;$L\geq L_{k. 0}$, then there exists \;$\lambda^L_{k-1}$\; such that for all \;;$\lambda\geq \lambda^L_{k-1}$, we have
$$
f_{k-1}(\lambda): B_{k-1}(\partial M) \longrightarrow  (\mathcal{E}_g)^{-L}
$$
is well defined and induces an isomorphism in homology.
\end{lem}
\vspace{6pt}

\noindent
Furthermore, still as in \cite{nd4}, we have also that Lemma \ref{eq:topnegative} and Lemma \ref{eq:mappingbary} imply the following one:
\begin{lem}\label{eq:verynegativesphere}
Assuming that  \;$k\geq 2$, $\bar k\geq 1$, $L\geq L_{k,,\bar k}$, then there exists \;$\lambda^L_{k-1, \bar k}>0$\; and \;$\Theta_{k-1, \bar k}^L>0$\; such that for all \;$\lambda\geq \lambda^L_{k-1, \bar k}$, we have
$$
\Psi_{k-1, \bar k}(\lambda, \Theta_{k-1, \bar k}^L): A_{k-1, \bar k} \longrightarrow  (\mathcal{E}_g)^{-L}
$$
is well defined and induces an isomorphism in homology.
\end{lem}

\subsection{Morse theoretical type results}
\begin{pfn} {of Theorem \ref{eq:morsepoincare1}-Theorem \ref{eq:Cm}}\\
The proof is the same as the one of Theorem 1.1-Theorem 1.6 in \cite{no1} by using Lemma \ref{eq:deformlemr}, Proposition \ref{eq:conspseudograd}, Corollary \ref{eq:loccritinf}, Lemma \ref{eq:morleminf}, Lemma \ref{eq:energybddinf}, Corollary \ref{eq:energybddcrit}, Lemma \ref{eq:tophigh}, and Lemma \ref{eq:topnegative} combined with the works of Bahri-Rabinowitz\cite{BR}, Karell-Karoui\cite{kk} and Malchiodi\cite{maldeg}. 
\end{pfn}
\subsection{Algebraic topological type results}

In order to carry the algebraic topological argument for existence, as in \cite{nd4}, we need the following lemma.
\begin{lem}\label{eq:nontriafee}
Assuming that $(ND)$ holds,  $s_k^*(O^*_{\partial M})\neq 0$\; in\; $H^3(S^{\infty})$\; and \;$s_k^*(O^*_{\partial M})=0$\; in\; $H^3(S_+^{\infty}\cup S_-^{\infty})$, then there exists $0\neq \tilde O^*_{\partial M}\in H^3(S)$ such that
 $$
i^*(\tilde O^*_{\partial M})= s_{k}^*(O^*_{\partial M}),
 $$
 where $i: S^{\infty}\longrightarrow S$ is the canonical injection.
\end{lem}
\begin{pf}
It follows from the same arguments as in the proof of Lemma 3.6 in \cite{nd4} by using the analysis of Section \ref{eq:cpi}.
\end{pf}

\begin{pfn}{ of Theorem \ref{eq:existence}}\\
The proof is the same as the one Theorem in \cite{nd4} by using the algebraic topological tools \eqref{eq:barytop}-\eqref{eq:defbarynega}, characterization of the critical points at infinity of \;$\mathcal{E}_g$\; established in Section \ref{eq:cpi}, and Lemma \ref{eq:nontriafee}.
\end{pfn}


\section{Appendix}\label{eq:appendix}

\begin{lem}\label{eq:intbubbleest}
Assuming that \;$\epsilon$\; is positive and small,  \;$a\in \partial M$\; and \;$\l\geq \frac{1}{\epsilon}$, then\\
1)
$$
\varphi_{a, \l}(\cdot)=\hat \d_{a, \l}(\cdot)+\log\frac{\l}{2}+H(a, \cdot)+\frac{1}{2\l^2}\D_{\hat g_{a}}H(a, \cdot)+O\left(\frac{1}{\l^3}\right)\;\;\;\text{on}\;\;\;\partial M
$$
2)
$$
\l\frac{\partial \varphi_{a, \l}(\cdot)}{\partial \l}=\frac{2}{1+\l^2\chi_{\varrho}^2(d_{\hat g_{a}}(a, \cdot))}-\frac{1}{\l^2}\D_{\hat g_{a}}H(a, \cdot)+O\left(\frac{1}{\l^3}\right)\;\;\;\text{on}\;\;\;\partial M,
$$
3)
$$
\frac{1}{\l}\frac{\partial \varphi_{a, \l}(\cdot)}{\partial a}=\frac{\chi_{\varrho}(d_{\hat g_{a}}(a, \cdot))\chi_{\varrho}^{'}((d_{\hat g_{a}}(a, \cdot))}{d_{\hat g_{a}}(a, \cdot)}\frac{2\l exp_a^{-1}(\cdot)}{1+\l^2\chi_{\varrho}^2(d_{\hat g_a}(a, \cdot))}+\frac{1}{\l}\frac{\partial  H(a, \cdot)}{\partial a}+O\left(\frac{1}{\l^3}\right);\;\; \text{on}\;\;\;\partial M,
$$
where \;$O(1)$ means $O_{a, \l, \epsilon}(1)$ and for it meaning see Section \ref{eq:notpre}.
\end{lem}
\begin{lem}\label{eq:outbubbleest}
Assuming that \;$\epsilon$\; is small and d positive, \;$a\in M$, $\l\geq \frac{1}{\epsilon}$, and \;$0<2\eta<\varrho$\; with \;$\varrho$\; as in \eqref{eq:cutoff}, then  there holds
$$
\varphi_{a, \l}(\cdot)=G(a, \cdot)+\frac{1}{2\l^2}\D_{\hat g_{a}}G(a, \cdot)+O\left(\frac{1}{\l^3}\right)\;\;\text{on} \;\;\;\partial M\setminus B^{a}_{a}(\eta),
$$
$$
\l\frac{\partial \varphi_{a, \l}(\cdot)}{\partial \l}=-\frac{1}{\l^2}\D_{\hat g_{a}}G_(a, \cdot)+O\left(\frac{1}{\l^3}\right)\;\;\text{on} \;\;\;\partial M\setminus B^{a}_{a}(\eta),
$$
and
$$
\frac{1}{\l}\frac{\partial \varphi_{a, \l}(\cdot)}{\partial a}=\frac{1}{\l}\frac{\partial G(a, \cdot)}{\partial a}+O\left(\frac{1}{\l^3}\right)\;\;\text{on} \;\;\;\partial M\setminus B^{a}_{a}(\eta),
$$
where \;$O(1)$\; means \;$O_{a, \l, \epsilon}(1)$\; and for it meaning see Section \ref{eq:notpre}.
\end{lem}
\vspace{4pt}

\noindent

\begin{lem}\label{eq:selfintest}
Assuming that \;$\epsilon$\; is small and positive, $a\in \partial M$ and $\l\geq \frac{1}{\epsilon}$, then there holds
$$
\mathbb{P}_g^{4, 3}\left(\varphi_{a, \l}, \;\varphi_{a, \l}\right)=16\pi^2\log \l-8\pi^2C_0+8\pi^2 H(a, a)+\frac{8\pi^2}{\l^2}\D_{\hat g_a} H(a, a)+O\left(\frac{1}{\l^3}\right),
$$
$$
\mathbb{P}_g^{4, 3}\left(\varphi_{a, \l},\; \l\frac{\varphi_{a, \l}}{\partial \l} \right)=8\pi^2-\frac{8\pi^2}{\l^2}\D_{\hat g_a}H(a, a)+O\left(\frac{1}{\l^3}\right),
$$
$$
\mathbb{P}_g^{4, 3}\left(\varphi_{a, \l}, \;\frac{1}{\l}\frac{\varphi_{a, \l}}{\partial a}\right)=\frac{8\pi^2}{\l}\frac{\partial H(a, a)}{\partial a}+O\left(\frac{1}{\l^3}\right),
$$
where \;$C_0$\; is a positive constant depending only on $n$,\;$O(1)$\; means \;$O_{a, \l, \epsilon}(1)$\; and for its meaning see Section \ref{eq:notpre}.
\end{lem}

\begin{lem}\label{eq:interactest}
Assuming that $\epsilon$ is small and positive $a_i, a_j\in \partial M$, \;$d_{\hat g}(a_i, a_j)\geq 4\ov C\eta$, $0<2\eta<\varrho$, $\frac{1}{\L}\leq \frac{\l_i}{\l_j}\leq \L$, and $\l_i, \l_j\geq \frac{1}{\epsilon}$, $\ov C$ as in \eqref{eq:proua}, and $\varrho$ as in \eqref{eq:cutoff}, then there hold
\begin{equation*}
\begin{split}
\mathbb{P}_g^{4, 3}\left(\varphi_{a_i, \l_i}, \;\varphi_{a_j, \l_j}\right)=&8\pi^2 G (a_j, a_i)+\frac{4\pi^2}{\l_i^2}\D_{\hat g_{a_i}} G(a_i, a_j)+\frac{4\pi^2}{\l_j^2}\D_{\hat g_{a_j}} G(a_j, a_i)\\&+O\left(\frac{1}{\l^3_i}+\frac{1}{\l_j^3}\right),
\end{split}
\end{equation*}
$$
\mathbb{P}_g^{4, 3}\left(\varphi_{a_i, \l_i}, \;\l_j\frac{\partial \varphi_{a_j, \l_j}}{\partial \l_j}\right)=-\frac{8\pi^2}{\l_j^2}\D_{\hat g_{a_j}}G(a_j, a_i)+O\left(\frac{1}{\l^3_j}\right),
$$
and
$$
\mathbb{P}_g^{4, 3}\left(\varphi_{a_i, \l_i}, \;\frac{1}{\l_j}\frac{\partial \varphi_{a_j, \l_j}}{\partial a_j}\right)= \frac{8\pi^2}{\l_j}\frac{\partial G(a_j, a_i)}{\partial a_j}+O\left(\frac{1}{\l^3_j}\right),
$$
where \;$O(1)$\; means here \;$O_{A, \bar \l, \epsilon}(1)$\; with \;$A=(a_i, a_j)$\; and \;$\bar \l=(\l_i, \l_j)$\; and for the meaning of $O_{A, \bar \l, \epsilon}(1)$, see Section \ref{eq:notpre}.
\end{lem}
\begin{lem}\label{eq:auxibubbleest1}
Assuming that $\epsilon>0$ is very small, we have that for $a\in \partial M$, $\l\geq \frac{1}{\epsilon}$, there holds
\begin{equation}
||\l\frac{\partial \varphi_{a, \l}}{\partial \l}||_{\mathbb{P}^{4, 3}}=\tilde O(1),
\end{equation}
\begin{equation}
||\frac{1}{\l}\frac{\partial \varphi_{a, \l}}{\partial a}||_{\mathbb{P}^{4, 3}}=\tilde O(1),
\end{equation}
and
\begin{equation}
||\frac{1}{\sqrt{\log \l}}\varphi_{a, \l}||_{\mathbb{P}^{4, 3}}=\tilde O(1),
\end{equation}
where here $\tilde O(1)$ means bounded by positive constants form below and above independent of $\epsilon$, $a$, and $\l$.
\end{lem}
\begin{lem}\label{l:A1}
1) If $\epsilon$ is small and positive, $a \in \partial M$, $p\in \N^*$, and $\l\geq \frac{1}{\epsilon}$ , then there holds
\begin{equation}
C^{-1}\l^{6p-3}\leq \oint_{\partial M}e^{3p\varphi_{a, \l}}dS_g\leq C\l^{6p-3},
\end{equation}
where $C$ is independent of $a$, $\l$, and $\epsilon$.\\\\
2) If $\epsilon$ is positive and small, $a_i, a_j \in \partial M$, $\l\geq \frac{1}{\epsilon}$ and $\l d_{\hat g}(a_i,a_j) \geq 4\ov C R$, then we have 
\begin{equation}
\mathbb{P}_g^{4, 3} (\varphi_{a_i,\l}, \;\varphi_{a_j,\l})\, \leq \,8\pi^2 G(a_i,a_j) \, + \, O(1),
\end{equation}
where \;$O(1)$\; means here \;$O_{A,  \l, \epsilon}(1)$\; with \;$A=(a_i, a_j)$, and for the meaning of $O_{A, \l, \epsilon}(1)$, see Section \ref{eq:notpre}. 
3) If $\epsilon$ is positive and small, $a_i, a_j \in \partial M$, $\l_i, \l_j\geq \frac{1}{\epsilon}$, $\frac{1}{\L}\leq\frac{\L_i}{\l_j}\leq \L$ and $\l _id_{\hat g}(a_i,a_j) \geq 4\ov C R$, then we have 
\begin{equation}
\mathbb{P}_g^{4, 3} (\varphi_{a_i,\l}, \;\varphi_{a_j,\l}) \, \leq \, 8\pi^2 G(a_i,a_j) \, + \, O(1),
\end{equation}
where \;$O(1)$\; means here \;$O_{A, \bar \l, \epsilon}(1)$\; with \;$A=(a_i, a_j)$\; and \;$\bar \l=(\l_i, \l_j)$\; and for the meaning of $O_{A, \bar \l, \epsilon}(1)$, see Section \ref{eq:notpre}. 
\end{lem}
\begin{lem}\label{l:A2}
Let  $p\in \N^*$,  $\hat R$ be a large positive constant, $\epsilon$ be a small positive number, $\alpha_i\geq 0$, $i=1, \cdots, p$, $\sum_{i=1}^p\alpha_i=k$, $\l\geq \frac{1}{\epsilon}$ and $ u = \sum_{i=1}^p \a_i \varphi_{a_i,\l}$. Assuming that there exist two positive integer $i, j\in \{1, \cdots, p\}$ with  $i\neq j$ such that $\l d_{\hat g}(a_i,a_j) \leq \frac{\hat R}{4\ov C}$, where $\ov C$ is as in \eqref{eq:proua}, then we have 
\begin{equation}
\mathcal{E}_g(u) \, \leq \mathcal{E}_g(v) \, + \, O(\log \hat R),
\end{equation}
with
$$
 v:= \sum_{k \leq p, k \neq i,j} \a_k \varphi_{a_k,\l} \, + (\a_i + \a_j) \varphi_{a_i,\l}.
$$
where here \;$O(1)$\;stand for \;$O_{{\bar\alpha}, A, \l, \epsilon}(1)$, with \;${\bar\alpha}=(\alpha_1, \cdots, \alpha_p)$\; and \;$A=(a_1, \cdots, a_p)$, and for the meaning of \;$O_{{\bar \alpha}, A, \l, \epsilon}(1)$, we refer the reader to Section \ref{eq:notpre}.
\end{lem}\begin{lem}\label{eq:lemma3}
1) If  $\epsilon$ is positive and small, $a_i, a_j\in \partial M$, $\l\geq \frac{1}{\epsilon}$ and $\l d_{\hat g}(a_i, a_j)\geq 4\ov C R$, then
$$
\varphi_{a_j, \l}(\cdot)=G(a_j,\cdot)+O(1)\;\;\;\text{in}\;\;B^{a_i}_{a_i}(\frac{R}{\l}),
$$
where here \;$O(1)$ means here \;$O_{ A, \l, \epsilon}(1)$, with  \;$A=(a_i, a_j)$, and for the meaning of \;$O_{ A, \l, \epsilon}(1)$, see Section \ref{eq:notpre}.\\
2) If $\epsilon$ is positive and small, $a_i, a_j\in \partial M$, $\l_i, \l_j\geq \frac{1}{\epsilon}$, $\frac{1}{\L}\leq\frac{\L_i}{\l_j}\leq \L$, and $\l_i d_{\hat g}(a_i, a_j)\geq 4\ov CR$, then
$$
\varphi_{a_j, \l_j}(\cdot)=G(a_j,\cdot)+O(1)\;\;\;\text{in}\;\;B^{a_i}_{a_i}(\frac{R}{\l_i}),
$$
where here \;$O(1)$ means here \;$O_{ A, \bar\l, \epsilon}(1)$, with  \;$A=(a_i, a_j)$, $\bar \l=(\l_i, \l_j)$ and for the meaning of \;$O_{ A, \l, \epsilon}(1)$, see Section \ref{eq:notpre}.\\
\end{lem}

\vspace{6pt}

\noindent
\begin{lem}\label{eq:positive}
There exists \;$\Gamma_0$\; and \;$\tilde \L_0$ two large positive constant such that for every $a\in \partial M$, $\l\geq \tilde \L_0$, and \;$w\in F_{a,  \l}:=\{w\in \mathcal{H}_{\frac{\partial }{\partial n}},  \ov{w}_{(Q, T)}=\left\langle \varphi_{a, \l}, w\right\rangle_{\mathbb{P}^{4, 3}}=\left\langle v_r, w\right\rangle_{\mathbb{P}^{4, 3}}=0,\;r=1\cdots,\bar k\}$, we have
\begin{equation}\label{eq:positivity}
 \oint_{\partial M}e^{3\hat\d_{a,\l}}w^2dV_{g_{a}}\leq \Gamma_0||w||^2_{\mathbb{P}^{4, 3}}.
\end{equation}
\end{lem}

\begin{lem}\label{eq:positiveg}
Assuming that\;$\eta$\; is a small positive real number with $0<2\eta<\varrho$ where $\varrho$ is as in \eqref{eq:cutoff}, then there exists a small positive constant $c_0:=c_0(\eta)$ and $\L_0:=\L_0(\eta)$ such that for every $a_i\in \partial M$ concentrations points  with $d_{\hat g}(a_i, a_j)\geq 4\ov C\eta$ where $\bar C$ is as in \eqref{eq:proua}, for every \;$\l_i>0$\; concentrations parameters satisfying $\l_i\geq \L_0$, with $i=1, \cdots, k$,  and for every \;$w\in E_{A, \bar \l}^*=\cap_{i=1}^k E^*_{a_i, \l_i}$ with $A:=(a_1, \cdots, a_k$), $\bar \l:= (\l_1, \cdots, \l_k)$ and $E^*_{a_i, \l_i}=\{w\in  \mathcal{H}_{\frac{\partial }{\partial n}}: \;\;\left\langle \varphi_{a_i, \l_i}, w\right\rangle_{\mathbb{P}^{4, 3}}=\left\langle \frac{\partial\varphi_{a_i, \l_i}}{\partial \l_i}, w\right\rangle_{\mathbb{P}^{4, 3}}=\left\langle \frac{\partial\varphi_{a_i, \l_i}}{\partial a_i}, w\right\rangle_{\mathbb{P}^{4, 3}}=\ov{w}_{(Q, T)}=\left\langle v_r, w\right\rangle_{\mathbb{P}^{4, 3}}=0,\;r=1\cdots,\bar k=0\}$, there holds
\begin{equation}\label{eq:positivity}
||w||^2_{\mathbb{P}^{4, 3}}- 6\sum_{i=1}^k\oint_{\partial M}e^{3\hat\d_{a_i,\l_i}}w^2dS_{g_{a_i}}\geq c_0||w||^2_{\mathbb{P}^{4, 3}}.
\end{equation}
\end{lem}


\begin{thebibliography}{99}
\bibitem{no1} Ahmedou M.,  Ndiaye C. B {\em Morse theory and the resonant \;$Q$-curvature problem}, Preprint: arxiv:.1409.7919.

\bibitem{akn} Ahmedou, M., Kallel S., Ndiaye, C. B., {\em The resonant boundary Q-curvature problem and boundary-weighted barycenters}. Preprint: arXiv:1604.03745.


\bibitem {aubin} Aubin T., Some nonlinear problems  in Riemannian geometry,  Springer Monographs in Mathematics, Springer-Verlag, Berlin 1998.\bibitem{ab1} Aubin T, Bahri A., {\em Methods of algebraic topology for the problem of prescribed scalar curvature}.  J. Math. Pures Appl. (9) 76 (1997), no. 6, 525-549.

\bibitem{ab2} Aubin T., Bahri A., {\em A topological hypothesis for the problem of prescribed scalar curvature}.  J. Math. Pures Appl. (9) 76 (1997), no. 10, 843-850.

\bibitem{bah} Bahri A., {\em Critical points at infinity in some variational problems}, Research Notes in Mathematics, 182, Longman-Pitman, London, 1989.

\bibitem{bah1} Bahri A., {\em Un probl\`eme variationel sans compacit\'e dans la g\'eom\'etrie de contact}, Comptes Rendus Math\'ematique Acad\'emie des Sciences, Paris, S\'erie I 299 (1984) 754-760. 

\bibitem{bah2}  Bahri A., {\em Pseudo-orbits of Contact Forms}, Pitman Research Notes in Mathematics Series, 173. Longman Scientific \& Technical, Harlow, 1988.

\bibitem{bd} Bahri A., {\em An invariant for Yamabe-type flows with applications to scalar-curvature problems in high dimension}. A celebration of John F. Nash, Jr. Duke Math. J. 81 (1996), no. 2, 323-466.

\bibitem{Bahri-Brezis} Bahri A., Brezis H., {\em Equations elliptiques non lin\'eaires sur des vari\'et\'es avec exposant de Sobolev critique}. C. R. Acad. Sci. Paris Sér. I Math. 307 (1988), no. 11, 573--576.

\bibitem{bc} Bahri A.,  Coron  J.M., {\em On a nonlinear elliptic equation involving the critical Sobolev exponent: the effect of the topology of the domain}, Comm. Pure Appl. Math. 41-3 (1988), 253-294.

\bibitem{BR} Bahri A., and Rabinowitz P., {\em Periodic solutions of Hamiltonian systems of $3$-body}, Ann. Inst. Poincar\'e, Anal. non lin\'eaire \textbf{8}(1991), 561-649.


\bibitem{bc} Bahri, A., Coron, J.M., {\em On a nonlinear elliptic equation involving the critical Sobolev exponent: the effect of the topology of the domain}, Comm. Pure Appl. Math. 41-3 (1988),-294.

\bibitem{bo} Branson T.P., Oersted., {\em Explicit functional determinants in four dimensions}, Proc. Amer. Math. Soc 113-3(1991), 669-682.

\bibitem{bred} Bredon G.E., {\em Topology and geometry},
Graduate Texts in Mathematics, 139, 1997.
cand., 57-2 (1995), 293-345.

\bibitem{bre1} Brendle S., {\em Curvature flows on surfaces with boundary}, Math. Ann. 324 (2002), no. 3, 491-519.

\bibitem{bre2} Brendle S., {\em A family of curvature flows on surfaces with boundary}, Math. Z. 241, 829-869 (2002).

\bibitem{bre3}  Brendle S., {\em A generalization of the Yamabe flow for manifolds with boundary}. Asian J. Math. 6 (2002), no. 4, 625–644.

\bibitem{bren3} Brendle S., {\em Convergence of the Yamabe flow for arbitrary initial energy}, J. Diff. Geom. 69 (2005),217-278.  

\bibitem{bs} Brendle S., Chen S. Z., {\em An existence theorem for the Yamabe problem on manifolds with boundary}. J. Eur. Math. Soc. (JEMS) 16 (2014), no. 5, 991-1016.

\bibitem{cn} Catino G., Ndiaye C. B., {\em Integral pinching results for manifolds with boundary}, Ann. Sc. Norm. Super. Pisa Cl. Sci. 9 (2010), 785-813.



\bibitem{cqy1} Chang S.Y.A., Qing J.,., Yang P.C.,{\em Compactification of a class of conformally flat $4$-manifold}, Invent. Math. 142-1(2000), 65-93.

\bibitem{cq1} Chang S.Y.A., Qing J.,{\em The Zeta Functional Determinants on manifolds with boundary 1. The Formula}, Journal of Functional Analysis 147, 327-362 (1997)

\bibitem{cq2} Chang S.Y.A., Qing J.,{\em The Zeta Functional Determinants on manifolds with boundary II. Extremal Metrics and Compactness of Isospectral Set}, Journal of Functional Analysis 147, 363-399 (1997)

\bibitem{cqy2} Chang S.Y.A., Qing J.,., Yang P.C.,{\em On the Chern-Gauss-Bonnet integral for conformal metrics on\;${\bf R}^{4}$},\;Duke Math. J.03-3(2000),523-544.

\bibitem{cl2} Chen C. C., Lin C. S, {\em Topological degree for a mean field equation on Riemann surfaces}, Comm. Pure, Appl. Math, 56-12 (2003), 1667-1727.
\bibitem{cs} Chen, S.S., {\em Conformal deformation on manifolds with boundary}, Geom. Funct. Anal. 19 (2009), no. 4, 1029--1064.


\bibitem{Dold} Dold A.,  {\em Algebraic topology}, Lectures on algebraic topology. Reprint of the 1972 edition. Classics in Mathematics. Springer-Verlag, Berlin, 1995. xii+377 pp. ISBN: 3-540-58660-1 55-02 (01A75).

\bibitem{dm} Djadli Z., Malchiodi A., {\em Existence of conformal metrics with constant \;$Q$-curvature}, Ann. of Math. (2) 168 (2008), no. 3, 813--858.


\bibitem{dr} Druet O., Robert F., {\em Bubbling phenomena for fourth-order four-dimensional PDEs with exponential growth}, Proc. Amer. Math. Soc 134(2006) no. 3, 897-908.

\bibitem{es1} Escobar J.F., {\em The Yamabe problem on manifolds with boundary}, Journal Differential Geometry 35 (1992) no.1, 21-84.









\bibitem{gt} Gilbar D., Trudinger N., {\em  Elliptic Partial Differential Equations of Second Order}, 2nd edition, Springr-Verlag, 1983.

\bibitem{gun} G\"unther M., {\em Conformal normal coordinates}, Ann. Global. Anal. Geom. 11 (1993), 173-184.



\bibitem{kk} Kallel S., Karoui R., {\em Symmetric joins and weighted barycenters}, Advanced Nonlinear Studies, 11(2011), 117--143.

\bibitem{lp}  Lee J., Parker T., {\em The Yamabe problem}, Bull. A.M.S. 17 (1987),37-81.

\bibitem{linwei} Lin C. S., Wei J., {\em Sharp estimates for bubbling solutions of a fourth order mean field equation}, Ann. Sc. Norm. Super. Pisa Cl. Sci. (5) 6 (2007), no. 4, 599-630.

\bibitem{lu} Lucia M., {\em A deformation lemma with an application to a mean field equation}. Topol. Methods Nonlinear Anal. 30 (2007), no. 1, 113--138.

\bibitem{lo} Lucia M., Or\'ak J., {\em A minimax theorem in the presence of unbounded Palais-Smale sequences}, Israel J. Math. 172 (2009).

\bibitem{maldeg} Malchiodi A., {\em Morse theory and a scalar field equation on compact surfaces}, Adv. Diff. Eq., 13 (2008), 1109-1129.

\bibitem{mal} Malchiodi A., {\em Compactness of solutions to some geometric fourth-order equations}, J. Reine Angew. Math. 594 (2006), 137--174.



\bibitem{marques} Marques F. C., {\em Existence results for the Yamabe problem on manifolds with boundary}, Indiana Univ. Math. J.  (2005),  1599-1620.

\bibitem{nd1} Ndiaye C. B., {\em Constant $Q$-curvature metrics in arbitrary dimension},  J. Funct. Anal. 251 (2007), no. 1, 1--58. 53.

\bibitem{nd2} Ndiaye C.B., {\em Conformal metrics with constant \;$Q$-curvature for manifolds with boundary}, Comm. Anal. Geom. 16 (2008), no. 5, 1049--1124.

\bibitem{nd6} Ndiaye C.B., {\em Constant \;$T$-curvature conformal metric on 4-manifolds with boundary},  Pacific J. Math. 240 (2009), no. 1, 151-184.

\bibitem{nd3} Ndiaye C. B., {\em $Q$-curvature flow on manifolds with boundary},  Math. Z. 269 (2011), 83-114. .

\bibitem{nd4} Ndiaye C. B., {\em Algebraic topological methods for the supercritical Q-curvature problem.} Adv. Math, 277, 277, pp 56--99, 2015.

\bibitem{nd7} Ndiaye C. B.,  {\em Topological methods for the resonant Q-curvature problem in arbitrary even dimension}, Journal of Geometry and Physics, Volume 140, June 2019, Pages 178-213.


\bibitem{nx} Ndiaye C. B., Xiao J, {\em An upper bound of the total $Q$-curvature and its isoperimetric deficit for higher-dimensional conformal Euclidean metrics}, Calc. Var. Partial Differential Equations, 38, no. 1-2, 1--27, 2010.


\bibitem{p1} Paneitz S., {\em A quartic conformally covariant differential operator for arbitrary pseudo-Riemannian manifolds}, preprint, 1983.

\bibitem{ops} Osgood B., Phillips R., Sarnak P., {\em Extremal of determinants of Laplacians}. J. Funct. Anal, 80, 148-211 (1988).


\bibitem{zw} Weinstein G., Zhang L., {\em The profile of bubbling solutions of a class of fourth order geometric equations on 4-manifolds}, J. Funct. Anal. 257 (2009), no. 12, 3895--3929

\end{thebibliography}
\end{document}